\documentclass[11pt]{article}
\usepackage[pctex32]{graphics}
\usepackage{amssymb}
\usepackage{hyperref}
\usepackage{latexsym}
\usepackage{epsfig} 
\usepackage{multirow} 
\usepackage{pstricks}
\usepackage{graphicx}
\usepackage{float}
\usepackage{epstopdf}
\usepackage{mathrsfs}
\usepackage{subcaption}
\usepackage{soul}
\usepackage{autonum}
\usepackage{amsmath}
\usepackage{amssymb}
\usepackage{relsize}
\usepackage{multirow}
\usepackage{color}
\usepackage{graphicx}
\usepackage{multicol} 

\usepackage{times}
\usepackage{epstopdf}
\newcommand{\comm}[1]{}
\newtheorem{definition}{Definition}[section]
\newtheorem{theorem}[definition]{Theorem}
\newtheorem{lemma}[definition]{Lemma}
\newtheorem{problem}[definition]{Problem}

\definecolor{Red}{rgb}{1,0.,0.}
\usepackage{fancyhdr}
\usepackage{setspace}

\newcommand{\R}{{\mathbb R}}

\newcommand{\mM}{{\mathsf M}}

\newcommand{\mK}{{\mathsf K}}

\newcommand{\mR}{{\mathsf R}}
\newcommand{\mT}{{\mathsf T}}
\newcommand{\mC}{{\mathsf C}}
\newcommand{\mI}{{\mathsf I}}
\newcommand{\mA}{{\mathsf A}}
\newcommand{\mV}{{\mathsf V}}
\newcommand{\mE}{{\mathsf E}}

\newcommand{\mL}{{\mathsf L}}
\newcommand{\mQ}{{\mathsf Q}}
\newcommand{\mO}{{\mathsf O}}
\newcommand{\mD}{{\mathsf D}}

\newcommand{\mU}{{\mathsf U}}

\newcommand{\mX}{{\mathsf X}}
\newcommand{\mH}{{\mathsf H}}
\newcommand{\mSigma}{{\mathsf \Sigma}}

\begin{document}

\title{Sparsity-promoting hierarchical Bayesian model for EIT with a blocky target}
\author{D Calvetti$^1$, M Pragliola$^2$, E Somersalo $^1$}
\date{$^1$ Case Western Reserve University, Cleveland, USA \\
$^2$ University of Naples Federico II, Naples, Italy}

\maketitle

\begin{abstract}
The electrical impedance tomography (EIT) problem of estimating the unknown conductivity distribution inside a domain from boundary current or voltage measurements requires the solution of a nonlinear inverse problem. Sparsity promoting hierarchical Bayesian models have been shown to be very effective in the recovery of almost piecewise constant solutions in linear inverse problems. We demonstrate that by exploiting linear algebraic considerations it is possible to organize the calculation for the Bayesian solution of the nonlinear EIT inverse problem via finite element methods with sparsity promoting priors in a computationally efficient manner. The proposed approach uses the Iterative Alternating Sequential (IAS) algorithm for the solution of the linearized problems. Within the IAS algorithm, a substantial reduction in computational complexity is attained by exploiting the low dimensionality of the data space and an adjoint formulation of the Tikhonov regularized solution that constitutes part of the iterative updating scheme. Numerical tests illustrate the computational efficiency of the proposed algorithm. The paper sheds light also on the convexity properties of the objective function of the maximum a posteriori (MAP) estimation problem.
\end{abstract}

\section{Introduction}

The goal in electrical impedance tomography (EIT) is to estimate the unknown electrical conductivity inside a body from current/voltage measurements at the boundary.  In the mathematical formulation of the problem, known as the Calder\'{o}n problem, the boundary data are expressed in terms of the Dirichlet-to-Neumann operator, mapping any applied voltage pattern to the corresponding current density pattern at the boundary. In practical applications, the measurements are performed by using a finite number of contact electrodes, and the data can be expressed in terms of the resistance matrix, mapping any applied current pattern vector to the corresponding voltage pattern. The connections between these two data types have been studied in the literature, see, e.g., \cite{hyvonen2009approximating,calvetti2019approximation}. For a relatively recent reviews and history of the EIT inverse problem, see, e.g., \cite{uhlmann2014inverse,adler2021electrical}.  In this article, the focus is on computational methods with electrode data.

It is not uncommon that a priori, there are reasons to believe that the unknown conductivity is nearly piecewise constant. Such belief can be imbedded in the reconstruction algorithms, e.g., by applying level set algorithms \cite{liu2017parametric}, using edge-enhancing penalties, including the total variation (TV) penalty or different variants of it \cite{borsic2012primal,harhanen2015edge,kaipio2000statistical,gonzalez2016experimental}. In the framework of Bayesian inverse problems, the edge-enhancing penalties can often be related to a sparsity-promoting prior, the conductivity being approximated by a discrete vector with sparse increments between values at adjacent loci.  In this work, we adopt this interpretation, considering a Bayesian hierarchical prior model that has been shown to yield a computationally efficient algorithm for edge enhancing in a number of linear inverse problems \cite{calvetti2019hierachical,calvetti2020sparse,calvetti2020sparsity}. 
The properties of the Iterative Alternating Sequential (IAS) algorithm have been analyzed extensively in the context of linear inverse problems in the cited articles, and has been recently applied to a non-linear inverse problem arising in diffuse optical tomography (DOT) in  \cite{manninen2022sparsity}, where further connections with other edge-enhancing methods have been discussed. 

One of the main questions in the present article is related to the computational efficiency. As pointed out in the previous articles, the IAS algorithm comprises two optimization tasks, the updating of the unknown of primary interest describing the conductivity, and the updating of the hyperparameter. The former problem requires a solution of a least squares problem that may be time consuming if the dimensionality of the problem is high. In previous works, a Krylov subspace method for fast approximation of the solution was proposed, and the viability was demonstrated by computed examples arising in linear inverse problems. In this article, the problem is revisited and it is shown that when the dimensionality of the data is significantly lower than the dimensionality of the unknown, the computational burden can be significantly reduced. The solution techniques proposed here are compared with earlier  approximate methods from the point of view of efficiency and accuracy. Moreover, we extend the theoretical discussion on the convergence of the IAS algorithm towards the non-linear model. In particular, we investigate to what extent the convergence analysis of the linear case remains valid for the type of non-linearity appearing in electrical impedance tomography. While the general convexity argument may not hold, the analysis leads to certain convexity results under appropriate limitations, shedding some light on the optimization landscape of the objective function. 

The article is organized as follows. In section~\ref{sec:model} the forward model is reviewed. While the finite element model is rather standard, we discuss some of the details of it as the structure of the stiffness matrix plays a role in the convexity analysis of the optimization problem. Sections~\ref{sec:bayesian} and \ref{sec:IAS} introduce and discuss Bayesian sparsity promoting hierarchical model and the maximum a posteriori estimation problem.  Section~\ref{sec:computational} discusses the numerical algorithms and  section~\ref{MAPvsqMAP} highlights the differences between finding or approximating the minimizer of the linearized objective function.  Finally, the convergence analysis is presented in section~\ref{sec:convergence}, and section~\ref{sec:computed} is dedicated to computed examples.

\section{Computational model}\label{sec:model}

We start by briefly reviewing the electrical impedance tomography model, as well the discretized version that constitutes the forward model, with the emphasis on the computational aspects. We then discuss the inverse problem from a similar perspective.

\subsection{Forward model}

Let $\Omega\subset\R^2$ be a bounded set with a connected complement, and let $\sigma:\Omega \to \R$ denote the electrical conductivity distribution. We assume that $0<\sigma_m\leq \sigma(x)\leq  \sigma_M<\infty$ for some positive constants $\sigma_m$ and $\sigma_M$.  We assume further that $L$ contact electrodes are attached to the boundary $\partial\Omega$, modeled as non-overlapping intervals $E_\ell\subset\partial\Omega$ of the boundary curve, $1\leq\ell\leq L$. 
Assuming that electric currents $I_\ell$ are injected through the electrodes, satisfying the Kirchhoff condition
\begin{equation}\label{Kirchhoff}
 \sum_{\ell = 1}^L I_\ell = 0,
\end{equation} 
the induced static electric voltage potential $u:\Omega\to \R$ satisfies the elliptic equation
\begin{equation}\label{voltage eq}
\nabla\cdot\big(\sigma\nabla u\big) = 0 \mbox{ in $\Omega$,}
\end{equation}
with the boundary conditions
\begin{equation}\label{boundary 1}
 \int_{E_\ell} \sigma\frac{\partial u}{\partial n} dS = I_\ell, \quad  \sigma\frac{\partial u}{\partial n} \bigg|_{\partial\Omega\setminus \cup E_\ell} = 0.
\end{equation}
This Neumann type boundary condition is not sufficient to determine uniquely $u$. To complement it, we assume that each electrode has a characteristic contact impedance $z_\ell >0$, and we impose the additional condition
\begin{equation}\label{boundary 2}
 \left(u + z_\ell\sigma\frac{\partial u}{\partial n}\right)\bigg|_{{E}_\ell} = U_\ell,
\end{equation}
where the constants $U_\ell$ are the electrode voltages that satisfy a ground condition, chosen here as
\begin{equation}\label{ground} 
 \sum_{\ell = 1}^L U_\ell = 0.
\end{equation}  
The forward problem is to solve the pair $(u, U)$ when  $\sigma$, $z$ and $I$ are given, where $U$, $I$ and $z$ are vectors in $\R^L$ with components $U_\ell$, $I_\ell$ and $z_\ell$, respectively.  The forward model is referred to as the complete electrode model (CEM) \cite{somersalo1992existence}. For an extension of the model with certain computational and theoretical advantages, we refer to \cite{hyvonen2017smoothened,darde2022contact}.

It is well known \cite{somersalo1992existence} that the CEM has a unique (weak) solution $(u,U)\in H^1(\Omega)\times \R^L_0$, where $\R^L_0$ is the set of real $L$-vectors satisfying the condition (\ref{ground}), and moreover, the solution satisfies the weak form equation,
\begin{equation}\label{weak form}
{\mathscr B}_{\sigma,z}\big((u,U),(v,V)\big) = \sum_{\ell=1}^L I_\ell V_\ell \quad\mbox{for all $(v,V) \in H^1(\Omega)\times \R^L_0$,}
\end{equation}
where the coercive quadratic form is given by
\[
{\mathscr B}_{\sigma,z}\big((u,U),(v,V)\big) = \int_\Omega \sigma\nabla u\cdot \nabla v + \sum_{\ell =1}^L \frac 1{z_\ell}\int_{{E}_\ell}(u-U_\ell)(v-V_\ell) dS .
\]
The above expression gives a natural finite element formulation of the forward model. Let ${\mathscr T}_h = \big\{K_\nu\big\}_{\nu = 1}^{n_t}$ denote a triangular tessellation of $\Omega$, such that $\overline\Omega_h = \cup \overline K_\nu$ is a polygonal approximation of the domain, where $h>0$ is a symbolic mesh size parameter. In the following, we do not distinguish between the domains $\Omega$ and $\Omega_h$.
Further, we denote by  $\{\psi_j\}_{j=1}^{n_v}$ a corresponding piecewise linear Lagrange basis, where $n_v$ is the number of vertices $p_k$ in the mesh, $\psi_j(p_k) = \delta_{jk}$.  We write an approximation of the voltage potential $u$ in $\Omega$,
\[
 u(x) = \sum_{j=1}^{n_v} u_j\psi_j(x),
\]
and let $\{{\mathcal{E}}_\ell\}_{\ell=1}^{L-1}$ denote a basis of $\R^L_0$, such that
\[
 U = \sum_{\ell=1}^{L-1}\alpha_\ell {\mathcal{E}}_\ell, \quad I = \sum_{\ell=1}^{L-1}\beta_\ell {\mathcal{E}}_\ell.
\] 
 By letting 
\[
 \big(u,U\big)  = \sum_{j=1}^{n_v} u_j \underbrace{\big(\psi_j,0\big)}_{\overline\psi_j} + \sum_{\ell=1}^{L-1} \alpha_\ell \underbrace{\big(0, {\mathcal{E}}_\ell\big)}_{\overline\psi_{n_v+\ell}}  = \sum_{j=1}^{n_v} u_j \overline\psi_j + \sum_{\ell = 1}^{L-1} \alpha_\ell \overline\psi_{n_v+\ell},
\] 
the weak form equation (\ref{weak form}) can be written component-wise as
\[
  \sum_{j=1}^{n_v} {\mathscr B}_{\sigma,z}\big(\overline\psi_j,\overline\psi_k \big)  u_j +  \sum_{\ell=1}^{L-1} {\mathscr B}_{\sigma,z}\big(\overline\psi_{n_v+\ell},\overline\psi_k \big) \alpha_\ell = \sum_{\ell = 1}^{L-1}\langle \overline\psi_{n_v+\ell},\overline\psi_k\rangle \beta_\ell, \quad 1\leq k\leq n_v+L,
 \]
or, in matrix form as
\begin{equation}\label{matrix}
 \left[\begin{array}{cc} \mK^{11} &  \mK^{12} \\ \mK^{21} & \mK^{22}\end{array}\right]\left[\begin{array}{c} u\\  \alpha \end{array}\right] = 
  \left[\begin{array}{c} 0\\  \beta \end{array}\right] ,
\end{equation}
where
\begin{eqnarray*}
 K_{jk}^{11}    &=& {\mathscr B}_{\sigma,z}\big(\overline\psi_j,\overline\psi_k \big),\quad 1\leq j,k\leq n_v, \\
 K_{j\ell}^{12} &=& {\mathscr B}_{\sigma,z}\big(\overline\psi_j,\overline\psi_{N+\ell} \big) = K^{21}_{\ell j},\quad 1\leq j\leq n_v, \; 1\leq \ell \leq L-1,\\
 K_{\ell \ell'}^{11} &=& {\mathscr B}_{\sigma,z}\big(\overline\psi_{n_v+\ell},\overline\psi_{n_v+\ell'} \big),\quad 1\leq \ell,\ell'\leq L-1.
 \end{eqnarray*}  
The components of the voltage vectors are solved from (\ref{matrix}) in terms of the the Schur complement of the block $\mK^{11}$ of $\mK$, often denoted by $\mK/\mK_{11}$, as
\[
 \alpha = \big(\mK/\mK_{11}\big)^{-1} \beta =  \left( \mK^{22} - \mK^{21}\big( \mK^{11} \right)^{-1}  \mK^{12} \big)^{-1}  \beta.
\]
To write the equation in terms of the electrode voltages and currents,  we define the matrix
\[
 \mE = \left[\begin{array}{ccc} {\mathcal{E}}_1 & \cdots & {\mathcal{E}}_{L-1}\end{array}\right] \in\R^{L\times (L-1)},
\]
and write $\beta$ in terms of the current vector $I$ as
\[
 \beta  = \big(\mE^\mT \mE\big)^{-1} \mE^\mT I = \mE^\dagger I,
\]
implying that
\[
 U = \mE \alpha = \underbrace{\mE \big(\mK/\mK_{11}\big)^{-1} \mE^\dagger}_{=\mR_{\sigma,z}} I  = \mR_{\sigma,z} I,
\]
where $\mR_{\sigma,z}\in\R^{L\times L}$ is the resistance matrix . Observe that, from the definitions, the dependencies on $\sigma$ and $z$ can be written as  $\mK^{11} = \mK^{11}(\sigma)$,   $\mK^{12} = \mK^{12}(z)$, $\mK^{21} = \mK^{21}(z)$ and  $\mK^{22} = \mK^{22}(z)$, indicating the non-linear dependency of the voltage vector on the conductivity and contact impedance. The significance of these dependencies will be elaborated further in the following section.

\subsection{Inverse problem}

To define the inverse problem, we assume that a set of $L-1$ linearly independent current patterns spanning the space $\R_0^L$,   $\{I^\ell\}_{\ell=1}^{L-1}$, referred to as a frame, is applied on the boundary of $\Omega$, and the corresponding set of voltages $\{U^\ell\}_{\ell=1}^{L-1}$ is measured. Assuming additive noise, the observation model can be written as
\[
 U^\ell = \mR_{\sigma,z} I^\ell + \varepsilon^\ell, \quad 1\leq \ell\leq L-1,
\]
where $\varepsilon^\ell$ is the observation noise.  The EIT inverse problem is to estimate $\sigma$, and sometimes $z$, from the measured voltages $U^\ell$.
Observe that while from the point of view of linear algebra, any linearly independent frame of current patterns can be chosen, the presence of the observation noise changes the situation dramatically, and a careful analysis of the sensitivity is necessary \cite{isaacson1986distinguishability}.

 To define the inverse problem, we assume that the conductivity can be written as
\[
 \sigma(x) = \sigma_0 + \delta\sigma(x),
\]
where $\sigma_0>0$ is a constant. Assume that there is a compact set $D\subset\Omega$ such that ${\rm supp}(\delta\sigma)\subset D$.  In particular, $\sigma = \sigma_0$ in the neighborhood $\Omega\setminus D$ of the boundary $\partial\Omega$, and we assume that the value $\sigma_0$ as well as the contact impedances $z_\ell$ are known. Furthermore, we assume that the set $D$ is conformal with the tessellation ${\mathscr T}_h$, that is, $\{1,2,\ldots,n_t\} = I_1 \cup I_2$ with disjoint index sets $I_1$ and $I_2$ such that
\[
 D = \underset{\nu\in I_1}\bigcup \overline K_\nu.
\]
Without loss of generality, we may assume that $I_1 = \{1,2,\ldots, n\}$ for some $n<n_t$. 

To define a finite dimensional version of the inverse problem, consider the dependency of $\mR_{\sigma,z}$ on $\sigma$. 
 In particular, the matrix $\mK^{11}(\sigma)$ depends on the conductivity through element-wise integrals,
\[
 I^\nu_{jk} = \int_{K_\nu} \sigma \nabla\psi_j\cdot \nabla \psi_k dx, \quad j,k\in {\rm vert}(K_\nu),
\]
where ${\rm vert}(K_\nu)$ is the triplet of indices of the three vertices of the triangle $K_\nu$.
Assuming that the basis functions are piecewise linear, the weak derivatives are piecewise constant, that is,
\[
  I^\nu_{jk} =  \left( \nabla \psi_j \cdot \nabla \psi_k \right) \big|_{K_\nu} \int_{K_\nu} \sigma  dx.
 \] 
 Based on this observation, we define the degrees of freedom of the inverse problem by
 \begin{equation}\label{xi}
  \xi_\nu = \frac 1{|K_\nu|}
  \int_{K_\nu} \delta\sigma dx, \quad 1\leq \nu\leq n,
 \end{equation}
 and reformulate the inverse problem as follows: {\em Assuming that $\sigma\in L^\infty(\Omega)$ with $0<\sigma_m\leq \sigma(x)\leq \sigma_M$ for some $\sigma_M>\sigma_m>0$, estimate the vector $\xi\in\R^{n}$ with components satisfying (\ref{xi}).} 
 
 In this formulation the conductivity is not discretized, and the information about it is contained in the discretized forward model depending on the tessellation ${\mathscr T}_h$. Changing the discretization merely changes the degrees of freedom, leaving the conductivity  invariant. Naturally, the discretization of the forward model introduces a modeling error. We say that two conductivities $\sigma$ and $\sigma'$ are {\em  ${\mathscr T}_h$-equivalent}, 
 \[
  \sigma \sim \sigma' ,\quad \mbox{if $\int_{K_\nu} \sigma dx = \int_{K_\nu} \sigma' dx$ for all $K_\nu$.}
 \]  
 Hence, if the degrees of freedom $\xi_\nu$ of $\sigma$ are given, a ${\mathscr T}_h$-equivalent conductivity that has the same degrees of freedom can be written as
 \[
 \sigma(x)  \sim  \sigma_h(x) = \sigma_0 + \sum_{\nu=1}^{n} \xi_\nu \chi_\nu(x), 
 \]
 where  $\chi_\nu$ is the characteristic function of the element $K_\nu$. However, we need not require that the underlying conductivity is piecewise constant.
 
 For later reference, we define the matrices $\mK_\nu^{11} \in\R^{n_v\times n_v}$ by 
 \begin{equation}\label{K11 nu}
  \big(\mK_\nu^{11} \big)_{jk} =  \int_{K_\nu} \nabla \psi_j \cdot \nabla \psi_k  dx =  |K_\nu| \left( \nabla \psi_j \cdot \nabla \psi_k \right) \big|_{K_\nu}, \mbox{ if $ j,k\in {\rm vert}(K_\nu)$ and $=0$ otherwise.}
 \end{equation}
 The matrix $ \mK^{11}$ is symmetric positive definite with ${\rm rank}\big(\mK_\nu^{11}) = 3$, and we observe that
 \[
  \mK^{11} = \sigma_0 \underbrace{\left( \sum_{\nu=1}^{n_t} K_\nu^{11}\right)}_{=\mK^{11}_0} + \sum_{\nu=1}^n \xi_\nu \mK_\nu^{11} = \sigma_0\mK_0^{11} +  \sum_{j=1}^n \xi_\nu \mK_\nu^{11} .
  \]
This formula constitutes the standard way of assembling the stiffness matrix, and, in addition, it is fundamental for computing the linearizations of the model.  
 
 \subsection{Discrete increments}
 
 Let  $\epsilon_\ell = \{v_j,v_k\}$ denote an edge directed from the node $v_j$ to $v_k$, $1\leq \ell \leq n_e$, where $n_e$ is the number of edges in the discretization mesh.
 We only consider edges in the interior or on the boundary of the set $D$. Let $N$ be the number of such edges. Let the edge $\epsilon_\ell$  separate two elements $K_\nu$ and $K_\mu$ not both included in $\Omega\setminus D$, that is,
 \[
  {\epsilon_\ell} = \overline K_\nu \cap\overline K_\mu, \quad (K_\nu \cup K_\mu)\cap D \neq \emptyset,
  \]
 with arbitrary orientation.  We define
 \[
  \iota_D(K_\nu) = 1, \mbox{if $\overline K_\nu \subset D$, and   $\iota_D(K_\nu) = 0$ otherwise,}
 \] 
and a sparse matrix $\mL\in\R^{N \times n}$ such that the only nonzero entries in the 
$\ell$th row are 
 \[
  \mL_{\ell \nu} = \iota_D(K_\nu), \quad \mL_{\ell \mu}  = -\iota_D(K_\mu).
 \]
The matrix has at most two non-zero entries on each row, with no zero rows. It is easy to see that 
${\rm rank}(\mL) = n<N$: In fact, if $v\in{\mathcal N}(\mL) \subset \R^n$, then by construction of the matrix $\mL$, the entries of $v$ corresponding to elements with neighboring elements in $\Omega\setminus D$ must vanish. It follows by induction that all the entries corresponding to interior elements in $D$ must vanish, hence,  ${\mathcal N}(\mL) = \{0\}$.

Given the matrix $\mL$, and the degrees of freedom $\xi$ of a conductivity, define a new increment variable,
\[
 \zeta = \mL \xi \in \R^N.
\]
Conversely, if $\zeta \in{\mathcal R}(\mL)$, it follows from the fact that $\mL$ has full rank, that 
\begin{equation}\label{zeta to xi}
 \xi = \big(\mL^\mT \mL\big)^{-1}\mL^\mT\zeta = \mL^\dagger \zeta, \quad \mL\xi = \zeta,
\end{equation}
where $\mL^\dagger$ is the pseudoinverse of $\mL$. 

In this work, we consider the EIT problem with the prior assumption that there is a ${\mathscr T}_h$-equivalent conductivity that can be explained in terms of few non-vanishing increments across the element boundaries. The problem can be formulated in terms of matrix $\mL$ as follows.

\begin{problem}
Assume that $\sigma = \sigma_0 + \delta\sigma \in L^\infty(\Omega)$ with $\delta\sigma$ supported in a compact set $D\subset \Omega$, and let  ${\mathscr T}_h$ be a finite triangular tessellation of the computational domain so that $D$ is conformal with it.
Given a complete set of current patterns $\{I^\ell\}_{\ell=1}^{L-1}$  and the noisy measurements of the corresponding voltages $\{V^\ell\}_{\ell=1}^{L-1}$, estimate the integrals (\ref{xi}) with the constraint that the vector $\zeta = \mL\xi \in\R^N$ is compressible, i.e., for a given threshold $\delta>0$,
\[
 \|\zeta\|_{0,\delta}  =\#\{j \mid 1\leq j\leq k, \|\zeta_\ell \|>\delta\} \ll N,
\]
that is, most of the components of $\zeta$ are below the threshold value.
\end{problem} 

In the following section, the problem is recast in terms of a hierarchical Bayesian model promoting the sparsity of the increment vector.

\section{Bayesian sparsity-promoting model}\label{sec:bayesian}

We consider the inverse problem in the Bayesian framework, thus modelling all unknowns as random variables, with the randomness reflecting the uncertainty about their values.  To set up the likelihood model, we arrange the observed voltages into a single vector,
\[
 b = \left[\begin{array}{c} U^1 \\ \vdots\\ U^{L-1}\end{array}\right] \in \R^m, \quad m = L(L-1),
\] 
and write the observation model
\[
 b = \left[\begin{array}{c} \mR_{\sigma,z} I^1 \\ \vdots \\  \mR_{\sigma,z} I^{L-1}\end{array}\right] + \left[\begin{array}{c} \varepsilon^1 \\ \vdots\\ \varepsilon^{L-1}\end{array}\right]
  = F(\xi,\sigma_0,z) + \varepsilon,
 \] 
 where the computational forward model is parametrized by the  integrals of $\sigma$ over the triangles. Since our main goal is estimating the variables $\xi$, we neglect the dependency on the presumably known $\sigma_0$ and $z$.
Assuming that the noise is a realization of a zero mean Gaussian with covariance matrix $\mSigma$,
 $\varepsilon \sim{\mathcal N}(0,\mSigma)$, the likelihood density can be written as
 \[
  \pi_{b\mid \xi}(b\mid \xi) \propto {\rm exp}\left( -\frac 12 \big\|\mSigma^{-1/2}\big(b - F(\xi)\big)\big\|^2\right).
 \] 
 To define the prior for the conductivity, we use the representation (\ref{zeta to xi}) of $\xi$ in terms of the increments over element edges, formulating the prior as a conditionally  Gaussian prior model,
 \begin{equation}\label{prior1}
  \pi_{\zeta\mid\theta}(\zeta\mid\theta) \propto \frac{1}{(\theta_1\cdots\theta_N)^{1/2}}{\rm exp}\left( - \frac 12 \sum_{\ell=1}^N \frac{\zeta_j^2}{\theta_j}\right), \quad \zeta \in {\mathcal R}(\mL).
 \end{equation}
 In other words, the prior density is a restriction of the $N$-dimensional Gaussian density ${\mathcal N}(0,\mD_\theta)$ in the range of the matrix $\mL$, where
 \[
  \mD_\theta = {\rm diag}(\theta_1,\ldots,\theta_N).
 \]  
 Observe that since if $\theta_j$ is small, then a priori the component $\zeta_j$ is assumed to be small, the prior model (\ref{prior1}) promotes sparsity if only few of the prior variances are  of significant size. Following \cite{calvetti2019hierachical,calvetti2023bayesian}, this goal can be achieved by assuming that $\theta$ is a random variable with mutually independent components with fat-tailed distributions, such as generalized gamma distributions, i.e.,
 \begin{equation}\label{prior2}
  \pi_\theta(\theta) = \prod_{j=1}^N \pi_{\theta_j}(\theta_j), \quad \pi_{\theta_j}(\theta_j) = \frac{|r|}{\Gamma(\beta)\vartheta_j}\left(\frac{\theta_j}{\vartheta_j}\right)^{r\beta -1}
 {\rm exp}\left( - \left(\frac{\theta_j}{\vartheta_j}\right)^r\right),
 \end{equation}
where $r\neq 0$, $\vartheta_j>0$ is a scale parameter, and $\beta>0$ is a shape parameter, the roles of which will be revisited later on. Combining formulas (\ref{prior1}) and (\ref{prior2}), we obtain the prior for the pair $(\zeta,\theta)$ as
\begin{equation}\label{joint prior}
\pi_{\zeta,\theta}(\zeta,\theta) = \pi_{\zeta\mid\theta}(\zeta\mid\theta)\pi_\theta(\theta) \propto {\rm exp}\left( - \frac 12 \sum_{j=1}^N \frac{\zeta_j^2}{\theta_j} - \sum_{j=1}^N\left(\frac{\theta_j}{\vartheta_j}\right)^r +\left(r\beta - \frac 32\right)\sum_{j=1}^N \log\frac{\theta_j}{\vartheta_j}\right), \quad \zeta\in {\mathcal R}(\mL).
\end{equation}
It follows from Bayes' formula that the posterior density can be written as
\[
 \pi_{\zeta,\theta\mid b}(\zeta,\theta\mid b) = {\rm exp}\big( - G(\zeta,\theta)\big), \quad  \zeta\in {\mathcal R}(\mL),
\]
with the Gibbs energy given by the formula
\begin{equation}\label{Gibbs}
  G(\zeta,\theta) = \frac 12 \big\|\mSigma^{-1/2}\big(b - F(\mL^\dagger\zeta)\big)\big\|^2  + \frac 12 \sum_{j=1}^N \frac{\zeta_j^2}{\theta_j} + \sum_{j=1}^N\left(\frac{\theta_j}{\vartheta_j}\right)^r  - \left(r\beta - \frac 32\right)\sum_{j=1}^N \log\frac{\theta_j}{\vartheta_j}.
\end{equation}
One of the aims of this work is the design of a computationally efficient way of estimating the maximum a posteriori (MAP) estimate, which is tantamount to minimizing the function $G(\zeta,\theta)$ under the constraint $\zeta \in  {\mathcal R}(\mL)$.

\section{MAP estimate by alternating iterations}\label{sec:IAS}

To minimize the Gibbs energy (\ref{Gibbs}), we use the Iterative Alternating Sequential (IAS) algorithm discussed in detail, e.g., in \cite{calvetti2023bayesian}. The minimization proceeds by alternating the following two steps until the convergence criterion is satisfied:
\begin{enumerate}
\item Fixing $\theta$ to its current value, minimize the part of $G$ that depends only on $\zeta$,
\begin{equation}\label{nonlin LSQ}
 G_\theta(\zeta) =  \frac 12 \big\|\mSigma^{-1/2}\big(b - F(\mL^\dagger\zeta)\big)\big\|^2  + \frac 12 \sum_{\ell=1}^N \frac{\zeta_j^2}{\theta_j} , \quad \zeta \in  {\mathcal R}(\mL).
\end{equation}
\item Fixing $\zeta$ to its current value, minimize the part of $G$ that depends only on $\theta$,
\[ 
G_\zeta(\theta) =  \frac 12 \sum_{\ell=1}^N \frac{\zeta_j^2}{\theta_j} + \sum_{j=1}^N\left(\frac{\theta_j}{\vartheta_j}\right)^r  - \left(r\beta - \frac 32\right)\sum_{j=1}^N \log\frac{\theta_j}{\vartheta_j}.
\] 
\end{enumerate}
The repetition of the updating steps generates a sequence of solutions $(\zeta^\ell,\theta^\ell)$, $\ell = 0,1,\ldots$. The iterations terminate as soon as
\begin{equation}\label{stopping}
\Delta_{\rm rel}\theta^\ell = \frac{\|\theta^\ell - \theta^{\ell-1}\|}{\|\theta^{\ell-1}\|} < \delta,
\end{equation}
where $\delta>0$ is a given tolerance.

We observe that the first step where $\zeta$ is updated is simply a non-linear least squares problem with a Tikhonov-type distributed regularization, while in the second step the updating can be performed for every component $\theta_j$ separately. 

Before addressing the computational details of the optimization problem, we summarize some known properties of the algorithm and elucidate the roles of the various hyperparameters in the model. 

\subsection{Insights from the IAS algorithm for linear problems}

The proposed iterative optimization algorithm has been extensively analyzed for linear inverse problems, and while the results do not extend automatically to non-linear problems like the one considered here, they shed some light on how to choose the hyperparameters in the model. In the following, we consider the Gibbs energy for a linear problem $F(\mL^\dagger\zeta) = \mA\zeta$,
\begin{equation}\label{Gibbs lin}
  \widetilde G(\zeta,\theta) = \frac 12 \big\|b - \mA \zeta\big\|^2  + \frac 12 \sum_{j=1}^N \frac{\zeta_j^2}{\theta_j} + \sum_{j=1}^N\left(\frac{\theta_j}{\vartheta_j}\right)^r  - \left(r\beta - \frac 32\right)\sum_{j=1}^N \log\frac{\theta_j}{\vartheta_j},
\end{equation}
where it is assumed, without loss of generality, that the noise is Gaussian white noise. 

We start by considering the case  $r=1$, corresponding to a gamma hyperprior. We assume that
\[
 \eta = \beta - \frac 32 >0.
\] 
 In \cite{calvetti2015hierarchical,calvetti2019hierachical}, it was proved that for $r=1$, the Gibbs energy (\ref{Gibbs lin}) is globally convex, hence has a unique minimizer and that the IAS algorithm converges to the global minimum. Moreover, it was shown  that  when $\eta  \rightarrow 0^+$,  the unique minimizer, denoted by $\zeta_\eta$, converges to $\zeta_0$,
 \[
  \zeta_0 = {\rm argmin}\left\{ \frac 12 \big\|b - \mA \zeta\big\|^2   + \sqrt{2}\sum_{j=1}^N \frac{|\zeta_j|}{\sqrt{\vartheta_j}}\right\}.
 \] 
This asymptotic result highlights two important facts concerning the hyperparameters. The first is that the value of the parameter $\eta>0$ can be set to control the sparsity, leading to the weighted $\ell^1$-norm estimate for small $\eta$. The second is that the hyperparameters $\vartheta_j$ play the role of weights in the penalty term. Setting the parameters $\vartheta_j$ all equal is tantamount to assuming a priori that the variables $(\zeta_j,\theta_j)$ are  independent and identically distributed, which is a special instance of exchangeability, i.e., the prior is invariant under permutations of the components. Exchangeability  expresses the prior belief that all components of the unknown are equally able to explain the data, an egalitarian position that in some classes of inverse problems, may lead to unsatisfactory solutions. In particular, it is well known that if the data are significantly more sensitive to some of the components than others, exchangeable priors favor solutions in terms of the sensitive components. For instance, in many geophysical and biomedical applications, exchangeable priors favor solutions in terms of pixel or voxel values near the sensors, leading to reconstructions with superficial sources.  In the current EIT application, we expect that with equal weighting, the algorithm favors explaining the data with conductivity perturbations in the proximity of the electrodes.

One successful way to avoid such biasing has been to weigh the variables by the corresponding sensitivities of the data, arguably a problematic practice in the Bayesian context.
Recently, in \cite{calvetti2019brain}, it was shown that by introducing the concept of ``signal-to-noise (SNR) exchangeability", stating that all variables should have the same possibility to explain the data {\em assuming that the SNR is fixed,} the sensitivity weighting can indeed be justified. This leads to the parameter choice 
\[
\vartheta_j = \frac{C}{\|a^{(j)}\|^2}, \quad a^{(j)} = \mA e_j,
\]
where $e_j\in\R^N$ is the canonical unit vector, and $C>0$ is a constant scalar that was discussed in the cited articles. We remark that for severely ill-posed problems, the norms of some columns of $\mA$ are vanishingly small, and the above formula may lead to values of $\vartheta_j$ that allow an unphysically large variance, requiring that an upper bound, independent of the exchangeability argument, for the values of $\vartheta_j$ needs to be imposed.
We observe that 
\[
 \|a^{(j)}\|  = \left\|\frac{\partial}{\partial \zeta_j}(\mA \zeta)\right\|  = \mbox{sensitivity of the forward model to $\zeta_j$,}
\] 
giving a natural interpretation that can be extended to non-linear problems such as EIT.

When $r<1$, the global convexity of the objective function is lost: a systematic study of local convexity can be found in \cite{calvetti2020sparse}. One reason for considering generalized gamma distributions with $r<1$ is the fact that, locally, the  convergence towards a sparse solution may be faster than for the case $r=1$. Typically, the values $r<1$ lead to sparser solutions than $r=1$, however there is no guarantee of finding a global minimizer. These observations led to the introduction of hybrid schemes \cite{calvetti2020sparsity} combining the conditionally Gaussian prior with different hypermodels.  

Hybrid algorithms are based on the idea to first run the IAS algorithm with the gamma hyperprior corresponding to $r=1$ to near convergence towards the unique global minimizer, then switch to a greedier  scheme by choosing a generalized gamma hyperprior with $r<1$.  We point out that while changing the prior model is not standard in the Bayesian philosophy, one may think of the different hypermodels as alternative parametrized expressions of the same underlying prior belief about the sparsity of the solution. To make the switching conceptually and computationally consistent, we match the hyperparameter values in the two models using the following conditions introduced in \cite{calvetti2023computationally},
\begin{enumerate}
\item If $\zeta_j=0$, we require that the values for $\theta_j$ in the two models coincide to guarantee that the a priori variance of the background outside the support of $\zeta$ is consistently defined independently of the model. 
\item The marginal expected value for $\theta_j$ is equal using both models.
\end{enumerate}
  
Denoting by $(r_1,\beta_1,\vartheta_1)$ and $(r_2,\beta_2,\vartheta_2)$ the hyperparameter values for two models with $r_1=1$,
these compatibility conditions imply that
\begin{equation}\label{compatibility}
 \vartheta_1\left(\beta_1 - \frac{3}{2r_1}\right)^{1/r_1} =  \vartheta_2\left(\beta_2 - \frac{3}{2r_2}\right)^{1/r_2},
\end{equation}
and by recalling the expectation of generalized gamma distribution, 
\begin{equation}\label{compatibility2}
\vartheta_1\frac{\Gamma(\beta_1 + \frac 1{r_1})}{\Gamma(\beta_1)} = \vartheta_2\frac{\Gamma(\beta_2 + \frac 1{r_2})}{\Gamma(\beta_2)}.
\end{equation} 
Observe that assuming a finite expectation of the generalized gamma distribution, certain restrictions of the hyperparameter values apply,

\subsection{Hyperparameter selection}

The analysis of linear inverse problems provides a guidance on how to set the model parameters. We begin running the hybrid IAS scheme, selecting $r=1$. To promote sparsity, we select $\eta = \beta_1 - 3/2$ to be small, e.g., $\eta = 1\times 10^{-5}$ and we set the hyperparameter $\vartheta$ based on the sensitivity analysis of the linear case. The sensitivity of the forward model to parameters $\zeta_j$ at $\zeta=0$ is
\[
 \left\|\frac{\partial}{\partial\zeta_j} F(\mL^\dagger\zeta)\big|_{\zeta=0}\right\| =\left\|D_\xi F(0)\mL^\dagger e_j\right\|,
\]   
therefore we set
\begin{equation}\label{vartheta}
 \vartheta_j = \frac{C}{\left\|D_\xi F(0)\mL^\dagger e_j\right\|^2},
\end{equation}
the choice of the parameter $C$  being discussed in the context of computed examples.

To run the hybrid IAS, we set $r_1=1$, and solve the pair $(\beta_2,\vartheta_2)$ from the compatibility conditions (\ref{compatibility})-(\ref{compatibility2}). For particular values of $r$, e.g., $r=-1$ or $r = \pm 1/2$,  explicit solutions can be found by using the properties of gamma functions. Again, the parameter values are given with the computed examples.

\section{Computational details of the IAS optimization}\label{sec:computational}

We consider now the implementational details of the IAS algorithm outlined in section~\ref{sec:IAS}. We begin with the second step, namely the optimization with respect to $\theta$, which is more straightforward.

\subsection{Updating of the variance}

Consider the second step of the IAS algorithm: Given the current value of $\zeta$, update the variance $\theta$.
To minimize the function $G_\zeta(\theta)$, consider the expression component-wise. We start by introducing the dimensionless variables: For $j$ arbitrary, define
\begin{equation}\label{scaling}
 t_j = \frac{\zeta_j}{\sqrt{\vartheta_j}}, \quad \lambda_j = \frac{\theta_j}{\vartheta_j},
\end{equation}
and temporarily neglecting the subindices, consider the function expressing the dependency of $G_\zeta$ on the scaled variables,
\[
  g(\lambda) = \frac 12\frac{t^2}{\lambda}  + \lambda^r - \eta \log\lambda, \quad \eta = r\beta - \frac 12.
\]  
Since as $\lambda \to 0+,\infty$, we have that $g(\lambda)\to\infty$, a minimizer exists and must satisfy the first order optimality condition $g'(\lambda) = 0$, or
\begin{equation}\label{critical}
   -\frac 12\frac{t^2}{\lambda^2}  + r\frac{\lambda^{r}}{\lambda} -\frac{\eta}{\lambda} = 0.
\end{equation}
For selected values of $r$, this equation admits a closed form solution.  The special parameter values include $r=1$ (gamma hyperprior), $r=-1$ (inverse gamma hyperprior), leading to
\[
 \lambda = \left\{\begin{array}{ll} \frac 12\left(\eta + \sqrt{\eta + 2 t^2}\right), & r = 1, \\ \frac 1{2|\eta|}(t^2 + 2), & r = -1.\end{array}\right.
\] 
For general values of $r\neq 0$, we need to resort to numerical approximations. By substituting an implicit functional representation of the solution, $\lambda = \Phi(t)$, in (\ref{critical})  and using implicit differentiation, we find that the function $\Phi$ must satisfy the initial value problem
\[
 \frac{d\Phi}{dt} = \frac{t\Phi}{t^2 + r(r-1)\Phi^r + \eta\Phi}, \quad \Phi(0) = \left(\frac{\eta}{r}\right)^{1/r}.
\]
By arranging the  numbers $t_j = \zeta_j/\sqrt{\vartheta_j}$ in increasing order, finding the corresponding values of $\lambda_j$ over the intervals $[t_j,t_j+1]$ starting with $t_0=0$ and using any numerical ODE solver leads to a fast updating of $\lambda_j$s and therefore of the variances $\theta_j$s.

\subsection{Updating the increment vector}

We consider now the first step, updating the increments  of the integrals of the conductivity over the elements using the current value of the variance $\theta$. Since finding the minimizer of (\ref{nonlin LSQ}) requires the solution of a nonlinear least squares problem of large dimensions with subspace constraint, to make the updating step fast and effective, extra care should be taken at each step.

We start by reparametrizing the problem, introducing an auxiliary variable $\alpha\in\R^N$ by
\[
 \alpha = \mD_\theta^{-1/2}\zeta \in {\mathcal R}(\mL_\theta), \quad \mL_\theta = \mD_\theta^{-1/2} \mL.
\]
and use  (\ref{zeta to xi}) to write $\alpha$ in terms of the original variable $\xi\in\R^n$, 
\[
 \alpha =  \mL_\theta \xi.
\]
Expressing the pseudoinverse of $\mL_\theta$ in terms of the QR-factorization of the matrix $\mL_\theta$,
\[
 \mL_\theta = \mQ \mR = \left[\begin{array}{cc}\mQ_1 & \mQ_2 \end{array}\right]  \left[\begin{array}{c}\mR_1 \\ \mO \end{array}\right], \quad \mQ_1 \in\R^{N\times n}, \;
 \mQ_2 \in\R^{N\times(N- n)},
\]
and observing that since $\mL$ is of full rank, the upper triangular matrix $\mR_1\in\R^{n\times n}$ is invertible,  we have 
\[
 \xi = \mR_1^{-1} \mQ_1^\mT \alpha = \mL_\theta^\dagger \alpha,\quad \mQ_2^\mT\alpha = 0,
\] 
and the solution of the minimization problem is equivalent to finding $\alpha$ that minimizes
\[
 \widehat G_\theta(\alpha) =  \frac 12 \big\|\mSigma^{-1/2}\big(b - F(\mL_\theta^\dagger\alpha)\big)\big\|^2  + \frac 12\|\alpha\|^2 , \quad \mQ_2^\mT \alpha = 0.
\] 
Since the first term of $\widehat G_\theta(\alpha)$ does not depend on $\mQ_2^\mT \alpha$, the minimizer of $\widehat G_\theta(\alpha)$ automatically satisfies the condition $\mQ_2^\mT \alpha = 0$. Therefore, the optimization problem is an unconstrained non-linear least squares problem. We note that the computation of the QR factorization, which may be very time consuming for large scale problems, could be a bottleneck for the algorithm. A discussion how this can be bypassed will be presented later. 

 We begin the numerical solution of  the problem by linearizing the functional at the current value of $\alpha$, starting with $\alpha_c = 0$. Given  $\alpha_c$, we  write
 $\alpha =\alpha_c + \delta\alpha$, and approximate the non-linear function as
 \begin{eqnarray*}
 F(\mL_\theta^\dagger\alpha) &\approx& F(\mL_\theta^\dagger\alpha_c) + D_\xi F(\mL_\theta^\dagger\alpha_c)\mL_\theta^\dagger \delta\alpha \\
 &=& F(\mL_\theta^\dagger\alpha_c) - D_\xi F(\mL_\theta^\dagger\alpha_c)\mL_\theta^\dagger \alpha_c +  
  D_\xi F(\mL_\theta^\dagger\alpha_c)\mL_\theta^\dagger \alpha. 
\end{eqnarray*} 
Letting
\begin{eqnarray*}
 \mA &=& \mSigma^{-1/2}D_\xi F(\mL_\theta^\dagger\alpha_c)\mL_\theta^\dagger \in\R^{m\times N},  \\
 r &=& \mSigma^{-1/2}\big(b - F(\mL_\theta^\dagger\alpha_c) + D_\xi F(\mL_\theta^\dagger\alpha_c)\mL_\theta^\dagger \alpha_c \big) \in\R^m,
\end{eqnarray*} 
the approximate the objective function (\ref{nonlin LSQ}) is of the form 
\begin{equation}\label{lin LSQ}
\widehat G_\theta(\alpha) = \frac 12 \big(\| r - \mA \alpha\|^2 + \|\alpha\|^2\big),
\end{equation}
which is the functional corresponding to a standard Tikhonov-regularized least squares problem with regularization parameter equal to one. The minimizer satisfies the corresponding normal equations,
\begin{equation}\label{normal}
 \big(\mA^\mT \mA + \mI_N\big)\alpha = \mA^\mT r,
\end{equation}
which is an $N\times N$ linear system. Since in our case $m<N$, a computationally more compact way is to solving the adjoint problem in data space,
\begin{equation}\label{adjoint}
 \big(\mA \mA^\mT + \mI_m\big) z = r, \quad \alpha = \mA^\mT z,
\end{equation}
which amounts to solving an $m\times m$ linear system. Indeed, if $z$ is the unique solution of the adjoint problem, by multiplying both sides by $\mA^\mT$ yields
\[
 \mA^\mT\big(\mA \mA^\mT + \mI_m\big) z = \big(\mA^\mT\mA + \mI_N\big)\mA^\mT z = \mA^\mT r,
\]
hence $\mA^\mT z$ must coincide with the unique solution of (\ref{normal}). Finally, for reasons of computational efficiency and numerical stability, the linear system (\ref{adjoint}) may be solved employing the  Lanczos bidiagonalization algorithm, discussed briefly below.

\subsubsection{Lanczos bidiagonalization process}

Recently, there has been a renewed interest in applying the Lanczos bidiagonalization process \cite{lanczos1952solution, paige1982lsqr}, or Golub-Kahan process \cite{golub1965calculating,golub2013matrix}, to solve large scale linear inverse problems, in particular in connection with hybrid regularization methods, see, e.g., \cite{chung2017generalized, chung2021computational,gazzola2020krylov,chung2023hybrid} and references therein. Krylov subspace iterative methods for the approximate solution of large scale linear systems of equations have several advantages over direct solvers based on matrix factorizations: The memory requirements for iterative solvers can be significantly smaller than for direct methods, making them the method of choice when memory allocation is an issue such as in embedded systems or when the linear system is prohibitively large, and they are the only feasible option when the linear forward map is given in a matrix-free form, allowing only the computation of the matrix-vector products.  While the computed examples presented in this work could be solved by direct methods in many computational environments, for completeness we include a brief review of the bidiagonalization process, because the formulation of the problem is particularly well-suited for that class of iterative solvers.

To summarize the algorithm, write the linear system (\ref{adjoint}) as
\begin{equation}\label{adjoint2}
 \big(\mM + \mI_m\big) z = r, \quad\mM = \mA \mA^\mT,
\end{equation}
and introduce the nested family of Krylov subspaces,
\[
 {\mathcal K}_\ell(r,\mM) = {\rm span}\big\{b,\mM b,\ldots,\mM^{\ell-1} b\big\}, \quad \ell\leq m.
\]
The $\ell$th approximate solution of the system (\ref{adjoint2}) satisfies
\begin{equation}\label{get z}
 z^{(\ell)} = \underset{z\in  {\mathcal K}_\ell(r,\mM)}{\rm argmin}\big\|(\mM+\mI_m)z - r\big\|^2.
\end{equation}
Given the matrices $\mA$ and $\mA^\mT$, or more generally, their actions on vectors, the Lanczos bidiagonalization process, summarized in algorithmic form below, produces a sequence  of matrices $\big\{\mC_\ell,\mU_\ell,\mV_\ell\big\}_{\ell=1}^m$,
\[
 \mU_\ell = \left[\begin{array}{ccc} u_1 & \cdots & u_\ell\end{array}\right]\in\R^{N\times\ell}, \quad  \mV_\ell = \left[\begin{array}{ccc} v_1 & \cdots & v_\ell\end{array}\right]\in\R^{m\times\ell}, \quad
\mC_\ell = \left[ \begin{array}{ c c c c c} \rho_1 &  & &  &\\ \sigma_2 & \rho_2 & & & \\ & \ddots & \ddots & &  \\ & & \sigma_{\ell-1} & \rho_{\ell-1} \\ &&&\sigma_\ell & \rho_\ell
\end{array}\right],
\]
satisfying the identities
\begin{eqnarray}
\mA \mU_\ell &= &\mV_\ell \mC_\ell + \sigma_{\ell+1}v_{\ell+1} e_\ell^\mT \label{LanczosIdentity1}
\\
\mA^\mT \mV_\ell &= &\mU_\ell \mC_\ell^\mT, \label{LanczosIdentity2}
\end{eqnarray}
where $e_\ell\in\R^\ell$ is the $\ell$th canonical unit vector, the columns of $\mU_\ell$ and $\mV_\ell$ are mutually orthogonal, and
\[
  {\mathcal K}_\ell(r,\mM) = {\rm span}\big\{v_1,\ldots,v_\ell\big\}, \quad r = \|r\| v_1.
\]  
Obviously, if $\ell=m$, $\sigma_{\ell+1}v_{\ell+1}e_\ell^\mT = 0$.  

In particular, if we write the solution of (\ref{get z}) in the form 
\[
 z  = \mV_\ell y,\quad y\in\R^\ell,
\]
we observe that from the identities (\ref{LanczosIdentity1}) and (\ref{LanczosIdentity2}) it follows that
\[
 \big(\mM+\mI_m\big)\mV_\ell y = \mV_\ell\big(\mC_\ell\mC_\ell^\mT + \mI_\ell\big)y + \sigma_{\ell+1}v_{\ell+1} e_\ell^\mT y,
\]
and further, the orthogonality of the vectors $v_j$ imply that
\[
 \big\|(\mM+\mI_m)z - r\big\|^2 =\big\|(\mC_\ell\mC_\ell^\mT + \mI_\ell) y - \|r\|e_1\big\|^2 + \sigma_{\ell+1}^2 y_\ell^2.
\] 
Thus, the problem of approximating the solution by $z^{(\ell)}$ can by reduced to solving an $\ell\times\ell$ tridiagonal problem, 
\[
(\mC_\ell\mC_\ell^\mT + \mI_\ell) y^{(\ell)}  = \|r\|e_1,
\]
where the  approximation error can be monitored by observing the quantity $\sigma_{\ell+1}|y_\ell|$. Finally, we point out that when memory allocation is an issue, the matrices $\mU_\ell$ and $\mV_\ell$ need not to be stored, however, passing from the solution $y^{(\ell)}$ to $z^{(\ell)}$ requires  repetition of the process, a trade-off of memory and complexity.  

\bigskip
\hrule
\medskip

{\bf Lanczos bidiagonalization algorithm}

\medskip
\hrule
\medskip

 \begin{enumerate}
 \item {\bf Given} $\mA \in\R^{m\times N}, \quad 0\neq r\in\R^m,\; 1<\ell <\min(m,N);$
 \item {\bf Initialize} $\sigma_1= \|b\|;\; v_1= b/\sigma_1; \;  \tilde{u} =\mA^\mT v_1; \; \rho_1= \| \tilde{u}\|; \; u_1 = \tilde{u}/\rho_1$.
 \item {\bf Repeat} for $j = 1,2,\ldots,\ell$
 \begin{enumerate}
 \item $\tilde{v}_j = \mA u_{j-1} - \rho_{j-1}v_{j-1}, \quad \sigma_j = \|\tilde{v}_j, \quad v_j = \tilde{v}_j/\sigma_j$;
 \item $\tilde{u}_j = \mA^\mT v_j - \sigma_j u_{j-1}, \quad \sigma_j = \|\tilde{u}_j \|, \quad u_j = \tilde{u}_j/\sigma_j$; 
 \end{enumerate}
 \item {\bf Compute}  $\tilde{v}_{\ell+1} = \mA u_\ell - \rho_\ell v_\ell,$ 
 \quad $\sigma_{\ell+1}= \|\tilde{v}_{\ell+1}\|, $
 \quad $v_{\ell+1}= \tilde{v}_{\ell+1} /\sigma_{\ell+1}$.
 \end{enumerate}

\medskip
\hrule
\medskip

\section{MAP and quasi-MAP estimates}\label{MAPvsqMAP}

In previous articles discussing the IAS problem for large scale problems, an approximate approach for estimating the minimizer of the linearized objective function (\ref{lin LSQ}) was proposed. Instead of seeking the true minimizer, the proposed approach made use of the fact that a viable alternative for finding a Tikhonov regularized solution for a linear ill-posed problem is to use Krylov subspace methods with an early stopping rule. More precisely, consider the sequence of least squares problem
\begin{equation}\label{CGLS}
 \alpha^{(\ell)} = \underset{\alpha \in{\mathcal K}_\ell(\mA^\mT r,\mA^\mT\mA)}{\rm argmin}\big\{\| r-\mA\alpha\|^2\big\}, \quad \ell=0,1,2,\ldots,
\end{equation}
constituting the sequence of iterates of the Conjugate Gradient for Least Squares (CGLS) method for solving the linear system $\mA\alpha =r$ in the least squares sense. When the matrix $\mA$ is ill-conditioned, as is the case in the present case, the least squares problem is ill-posed, requiring regularization such as the Tikhonov regularization. A viable alternative is to stop the iterations for updating the approximations (\ref{CGLS}) before the noise takes over \cite{hanke2001lanczos}. Assuming that the noise $\varepsilon\in\R^m$ in the (whitened) data is a realization from the standard normal distribution, the stopping criterion can be to let $\ell$ be the smallest integer satisfying 
\[
  \|r - \mA\alpha^{(\ell)}\|^2 \leq m = {\mathbb E}(\|\varepsilon\|^2),
\]  
or, alternatively,
\[
 \widehat G(\alpha^{(\ell+1)}) > \widehat G(\alpha^{(\ell)}).
\]
The latter condition amounts to stopping the iterates at the minimum of the semiconvergence of the objective function, where further iterations would cause it to increase. 

\section{Analysis of the convergence}\label{sec:convergence}

As pointed out in the previous sections, the IAS algorithm with the hypermodel $r=1$ and linear forward model is globally convergent to a unique minimizer, which follows from the global convexity of the objective function to be minimized. In the case of a non-linear forward model, such conclusion does not necessarily hold, however, owing to the special structure of the forward problem, a partial result can be proved. Some of the results shown here are of interest also for the convergence of EIT optimization algorithms based on standard Tikhonov regularization.  

For the sake of definiteness, let us assume that the basis vectors $\{\mathcal{E}_\ell\}_{\ell=1}^{L-1}$ of $\R^L_0$ are mutually orthonormal, and the current patterns $\{I^\ell\}_{\ell=1}^{L-1}$ applied on the electrodes are $I^\ell = \mathcal{E}_\ell$. Therefore, writing $I^\ell = \mE \beta^\ell$, we have $\beta^\ell = e_\ell\in\R^{L-1}$, $1\leq \ell\leq L-1$, and
the forward model (\ref{matrix}) to be solved for a full frame of current patterns is
\[
 \mK \left[\begin{array}{cccc} u^1 & u^2 &\cdots & u^{L-1} \\ \alpha^1 & \alpha^2 &\cdots &\alpha^{L-1}\end{array}\right] = 
 \left[\begin{array}{cccc} 0 & 0 &\cdots & 0 \\ \beta^1 & \beta^2 &\cdots & \beta^{L-1} \end{array}\right]  = \left[\begin{array}{c} \mO_{N\times(L-1)} \\ \mI_{L-1}\end{array}\right],
\]
or, concisely,
\begin{equation}\label{mateq}
 \mK \mX = \left[\begin{array}{c} \mO \\ \mI_{L-1}\end{array}\right], \quad \mX =  \left[\begin{array}{cccc} u^1 & u^2 &\cdots & u^{L-1} \\ \alpha^1 & \alpha^2 &\cdots &\alpha^{L-1}\end{array}\right] .
\end{equation}
In the following, we consider both $\mK$ and $\mX$ as matrix valued functions of the parameters $\xi_\nu$, $1\leq\nu\leq n$ as defined in (\ref{xi}).  To analyze the convexity of the objective function of the
iterative algorithm, we need to calculate the derivatives up to the second order of $\alpha\in\R^{(L-1)\times(L-1)}$ with respect to these parameters.

Differentiating both sides of (\ref{mateq}) with respect to $\xi_\nu$ gives
\begin{equation}\label{leibniz}
 \big(\partial_\nu\mK\big) \mX + \mK  \big(\partial_\nu\mX\big) =0,
\end{equation}
from which it follows that
\[
 \partial_\nu \mX = - \mK^{-1}\big(\partial_\nu \mK\big) \mX.
\]
Multiplication from the left by the transpose of the right hand side in  (\ref{mateq}) yields
\[
 [\mO^\mT \; \mI_{L-1}]  \partial_\nu \mX =  - [\mO^\mT \; \mI_{L-1}]  \mK^{-1}\big(\partial_\nu \mK\big) \mX.
\] 
Using the symmetry of $\mK$, it follows that
\[
 [\mO^\mT \; \mI_{L-1}]  \mK^{-1}\big(\partial_\nu \mK\big) \mX = \left( \mK^{-1} \left[\begin{array}{c} \mO \\ \mI_{L-1}\end{array}\right]\right)^\mT\big(\partial_\nu \mK\big) \mX
  = \mX^\mT \big(\partial_\nu \mK\big) \mX,
\]
thus yielding the formula
\begin{eqnarray}
 \partial_\nu\alpha &=& \partial_\nu \left[\begin{array}{cccc} \alpha^1 &\alpha^2 &\cdots &\alpha^{L-1}\end{array}\right] =
  [\mO^\mT \; \mI_{L-1}]  \partial_\nu \mX \nonumber\\
  &=& -\mX^\mT \big(\partial_\nu \mK\big) \mX,
  \label{dalpha}
\end{eqnarray}
that is used to compute numerically the Jacobian of the forward map. Observe, in particular, that
\[
 \partial_\nu\mK = \left[\begin{array}{cc} \partial_\nu \mK^{11} & \mO \\ \mO & \mO\end{array}\right] =  \left[\begin{array}{cc} \mK^{11}_\nu & \mO \\ \mO & \mO\end{array}\right],
\]
where $\mK_\nu^{11}\in\R^{n\times n}$ is a symmetric rank-3 matrix (\ref{K11 nu}) independent of $\xi$.

To compute the second derivatives,  we differentiate both sides of equation (\ref{dalpha}) with respect to $\xi_\mu$ and use the fact that $\partial_\nu \mK$ is independent of $\xi$ to obtain
\[
 \partial_{\mu\nu}^2\alpha =  - (\partial_\mu \mX)^\mT (\partial_\nu\mK) \mX - \mX^\mT(\partial_\nu \mK) \partial_\mu \mX. 
\] 
 From (\ref{leibniz}) we have 
\[
 (\partial_\nu \mK) \mX = - \mK(\partial_\nu \mX),
\]
thus
\begin{equation}\label{d2alpha}
   \partial_{\mu\nu}^2\alpha = (\partial_\mu X)^\mT\mK(\partial_\nu X) + (\partial_\nu\mX)^\mT \mK (\partial_\mu\mX).
\end{equation} 
We consider now the data-dependent part of the objective function (\ref{Gibbs}). To simplify the notations, let us assume for the time being that the noise is scaled white noise, $\varepsilon\sim{\mathcal N}(0,\omega^2\mI)$. Furthermore, we write the measured voltage vectors in terms of the orthonormal basis $\mE$ as
\[
 U^\ell = \mE \gamma^\ell, \quad 1\leq \ell \leq L-1,
\] 
so that
\[
 \left[\begin{array}{ccc} U^1 & \cdots &U^{L-1}\end{array}\right] = \mE \left[\begin{array}{ccc} \gamma^1 & \cdots &\gamma^{L-1}\end{array}\right] = \mE\gamma.
\]  
With these notations, we write
\[
 g(\xi) = \frac 1{2\omega^2}\|b - F(\xi)\|^2 =  \frac 1{2\omega^2}\|[\mE\gamma - \mE\alpha(\xi)\|^2_F = \frac 1{2\omega^2}\|\gamma - \alpha(\xi)\|^2_F 
\]
by the orthogonality of the basis $\{{\mathcal E}_\ell\}_{\ell=1}^{L-1}$, $\|\,\cdot\,\|_F$ being the Frobenius norm.
To analyze the convexity properties of the objective function, we need to calculate the second derivatives of $g$. Differentiating, we have
\[
 \partial_\nu g(\xi) = \frac 1{2\omega^2}\sum_{i,j=1}^{L-1} (\alpha_{ij} - \gamma_{ij})\partial_\nu\alpha_{ij},
\] 
and
\begin{equation}\label{C and D}
\partial_{\mu\nu}^2 g(\xi) = \frac 1{2\omega^2}\underbrace{\sum_{i,j=1}^{L-1} \partial_\nu \alpha_{ij}\partial_\mu\alpha_{ij}}_{=D_{\mu\nu}} + \frac 1{2\omega^2}\underbrace{\sum_{i,j=1}^{L-1} (\alpha_{ij} - \gamma_{ij})\partial^2_{\mu\nu}\alpha_{ij}}_{=C_{\mu\nu}} = \frac{1}{2\omega^2}\left(D_{\mu\nu}+C_{\mu\nu}\right).
\end{equation}
Observe that while the matrix $\mD = \big[D_{\mu\nu}\big]$ does not depend on data, the matrix $\mC = \big[C_{\mu\nu}\big]$ does, and due to the presence of the noise in data, it is a random matrix. In the following lemma, we consider only the expectation of it, writing
\[
 \overline C_{\mu\nu} = {\mathbb E}( C_{\mu\nu}).
\] 

\begin{lemma}
The matrix $\mD =\left[D_{\mu\nu}\right]\in\R^{n\times n}$ is symmetric positive semidefinite. Moreover, if $\overline b = {\mathbb E}(b)$ is the noiseless voltage vector that corresponds to an underlying true conductivity
\begin{equation}\label{sigma star}
 \sigma^* = \sigma_0 + \sum_{\nu=1}^{n_t} \xi_\nu^*\chi_\nu, \mbox{ and $\xi_\nu\leq \xi_\nu^*$ for all $\nu$}, 
\end{equation}
then the matrix  $\overline \mC = \left[\overline C_{\mu\nu}\right]$ is symmetric positive semidefinite.
\end{lemma}

{\bf Proof:} Let $v\in\R^{n}$, $v\neq 0$.  Denoting by $\partial_v$ the derivative in the direction of $v$, we have
\[
 v^\mT \mD v = \sum_{\mu,\nu=1}^{n_t} \sum_{i,j =1}^{L-1} v_\nu \partial_\nu \alpha_{ij} v_\mu\partial_\mu\alpha_{ij} =  \sum_{i,j =1}^{L-1}  \partial_v \alpha_{ij} \partial_v\alpha_{ij} = \|\partial_v\alpha\|_F^2\geq 0,
\] 
proving the positive semidefiniteness of $\mD$. 

To prove the statement for the matrix $\overline\mC$, we first observe that by (\ref{d2alpha}),
\[
 \sum_{\mu,\nu=1}^{n} v_\mu v_\nu \partial_{\mu\nu}^2\alpha = 2 (\partial_v\mX)^\mT \mK(\partial_v X),
\] 
which is a symmetric positive semidefinite matrix in $\R^{(L-1)\times(L-1)}$.  Defining 
\begin{equation}\label{Delta}
\Delta = \alpha-\gamma,
\end{equation}
we see that
\[
 v^\mT\mC v = 2 \left\{(\partial_v\mX)^\mT \mK(\partial_v X)\right\}:\Delta = \left\{(\partial_v\mX)^\mT \mK(\partial_v X)\right\}:(\Delta + \Delta^\mT),
 \]
 where we used the notation $``: "$ for the Frobenius inner product for matrices in $\mR^{(L-1)\times (L-1)}$. 
 The symmetric positive semidefiniteness of the first matrix in the right hand side allows us to express it in terms of its eigenvalue decomposition as
 \[
  (\partial_v\mX)^\mT \mK(\partial_v X) = \sum_{\ell=1}^{L-1}\eta_\ell w_\ell w_\ell^\mT, \quad \eta_\ell\geq 0,\quad w_\ell\in\R^{L-1},
 \]
therefore
 \[
   v^\mT \mC v = \sum_{\ell=1}^{L-1} \eta_\ell w_\ell^\mT(\Delta + \Delta^\mT) w_\ell.
  \] 
  The expected value of $ v^\mT \mC v $ is guaranteed to be non-negative if
  \begin{equation}\label{exp C}
  {\mathbb E}\big( w^\mT(\Delta + \Delta^\mT) w\big)\geq 0
 \end{equation}
 for all $w\in\R^{L-1}$. 
 Consider first the matrix $\alpha$ in the expression (\ref{Delta}). Indicating explicitly the dependency of the matrix $\mK$ on the conductivity, we have
 \[
 \left[\begin{array}{c} u \\ \alpha\end{array} \right] =  \mK_\sigma^{-1} \left[\begin{array}{c} \mO \\ \mI_{L-1}\end{array}\right],
 \]
therefore   
\begin{eqnarray*}
 w^\mT \alpha w  &=&  \left[\begin{array}{cc} 0 & w^\mT\end{array}\right]  \left[\begin{array}{c} u \\ \alpha\end{array}\right] w  \\ 
 &=&   \left[\begin{array}{cc} 0 & w^\mT\end{array}\right] \mK_\sigma^{-1}  \left[\begin{array}{c} \mO \\ \mI_{L-1}\end{array}\right] w \\
  &=& \left[\begin{array}{cc} 0 & w^\mT\end{array}\right] \mK_\sigma^{-1}  \left[\begin{array}{c} 0 \\ w \end{array}\right].
\end{eqnarray*}

 A similar argument with the assumption that the noiseless data arise from an underlying conductivity $\sigma^*$ implies that the expectation $\overline\gamma = {\mathbb E}(\gamma)$ satisfies
\[
 w^\mT \overline\gamma w =  \left[\begin{array}{cc} 0 & w^\mT\end{array}\right] \mK_{\sigma^*}^{-1}\left[\begin{array}{c} 0 \\ w\end{array}\right].
\] 
Therefore, by the symmetry of the matrices $\mK_\sigma$ and $\mK_{\sigma^*}$,
\[
 {\mathbb E}\big(w^\mT(\Delta + \Delta^\mT) w \big) = 2 \left[\begin{array}{cc} 0 & w^\mT\end{array}\right] \big(\mK_\sigma^{-1} - \mK_{\sigma^*}^{-1}\big)\left[\begin{array}{c} 0 \\ w\end{array}\right] .
\] 
This expression can be interpreted to represent the projection into direction $w$ of the coefficient vector of the voltage difference between targets of conductivities $\sigma$ and $\sigma^*$ when applying a current pattern 
$I = \mE w$  at the electrodes.
It follows from the  assumption $\sigma\leq \sigma^*$, that
\[
 K_\sigma \leq K_{\sigma^*}
\]
in the sense of quadratic forms, hence
 \[
 K_{\sigma^*}^{-1} \leq K_{\sigma}^{-1},
\]
showing that the inequality (\ref{exp C}) hold, thus completing the proof.
$\Box$

Consider the the Hessian of the Gibbs energy functional (\ref{Gibbs}) with respect to the pair $(\zeta,\theta)$ in the case where $r=1$. We  partition it as  
\[
 \mH = \left[\begin{array}{cc} \nabla_\zeta\nabla_\zeta G(\zeta,\theta) & \nabla_\zeta\nabla_\theta G(\zeta,\theta) \\ \nabla_\theta\nabla_\zeta G(\zeta,\theta) & \nabla_\theta\nabla_\theta G(\zeta,\theta)\end{array}\right]
\]
and analyze the structure of each block.
A straightforward differentiation shows that
\[
 \nabla_\theta \nabla_\theta G(\zeta,\theta) = {\rm diag}\left(\frac{\zeta_j^2}{\theta_j^3} + \frac{\eta}{\theta_j^2}\right),
\] 
and
\[
 \nabla_\zeta\nabla_\theta G(\zeta,\theta) = {\rm diag}\left( -\frac{\zeta_j}{\theta_j^2}\right).
\]
Assuming for simplicity scaled white noise, 
\begin{eqnarray*}
    \nabla_\zeta\nabla_\zeta G(\zeta,\theta) &=&  \nabla_\zeta\nabla_\zeta\left(\frac 1{2\omega^2}\|b - F(\mL^\dagger\zeta)\|^2 + \frac 12 \sum_{j=1}^N\frac{\zeta_j^2}{\theta_j}\right) \\
    &=& \nabla_\zeta\nabla_\zeta g(\mL^\dagger \zeta) + {\rm diag}\left(\frac 1{\theta_j}\right) \\
    &=& \frac 1{2\omega^2}\left(\mL^\dagger\right)^\mT(\mC + \mD)\mL^\dagger + {\rm diag}\left(\frac 1{\theta_j}\right),
 \end{eqnarray*}   
 where the matrices $\mC$ and $\mD$ are defined in (\ref{C and D}).

To test the positive definiteness of the Hessian, let $q = [q^1;q^2]\in\R^{2N}$, $q\neq 0$. Then 
\begin{eqnarray*}
q^\mT \mH q &=&  \left(\mL^\dagger q^1\right)^\mT(\mC + \mD)\mL^\dagger q^1 + \sum_{j=1}^N \frac{(q_j^1)^2}{\theta_j}  - 2\sum_{j=1}^N\frac{\zeta_j}{\theta_j^2}{q^1_jq^2_q}  + \sum_{j=1}^N\left(\frac{\zeta_j^2}{\theta_j^3} + \frac{\eta}{\theta_j^2}\right) (q^2_j)^2 \\
&=& \frac 1{2\omega^2} \left(\mL^\dagger q^1\right)^\mT(\mC + \mD)\mL^\dagger q^1 +\sum_{j=1}^2\frac 1{\theta_j}\left(q_j^1 -\frac{\zeta_jq_j^2}{\theta_j}\right)^2 + 
 \sum_{j=1}^N \frac{\eta}{\theta_j^2} (q^2_j)^2,
\end{eqnarray*}
 showing that under the assumption of the previous lemma, the expectation of the Hessian is symmetric positive definite. Observe that in the case of a linear forward model, the matrix $\mC$ would vanish, leading to a global convexity of the objective function. This observation was the basis of the proof of uniqueness of the global minimizer as well as convergence of the IAS algorithm in the linear case. In the current problem, such conclusion is not possible, however, a partial convexity result is summarized in the following theorem.

\begin{theorem}
Given a noiseless voltage measurement that corresponds to a conductivity $\sigma^* =\sigma_0 + \sum_{\nu=1}^n\xi_\nu^*\chi_\nu$, a sufficient conditions to guarantee that the Gibbs energy $G(\zeta,\theta)$ with $r=1$ is convex is  $\xi_\nu \leq \xi_\nu^*$, $1\leq \nu \leq n$, where $\xi = \mL^\dagger \zeta$.
\end{theorem}

We remark that the above theorem gives only a sufficient condition for the convexity. In general, the analysis in this section allows the explicit computation of the Hessian and a numerical verification of the convexity. 

\section{Computed examples}\label{sec:computed}

In all computed examples, the domain representing the body $\Omega$ is the unit disc. A set of 32 identical electrodes are attached on the boundary of $\Omega$, distributed uniformly around the unit circle. The space filling fraction, defined as the fraction of the boundary covered by the electrodes is 0.45. Furthermore, we assume that the conductivity has a known background value $\sigma_0=1$ outside the compact domain $D$, defined as the closed disc centered at the origin and having a radius $R=0.9$. Hence, the inverse problem is to estimate the conductivity in the disc $D$. 

\begin{figure}[ht!]
\centerline{
\includegraphics[width=6.2cm,angle=-90]{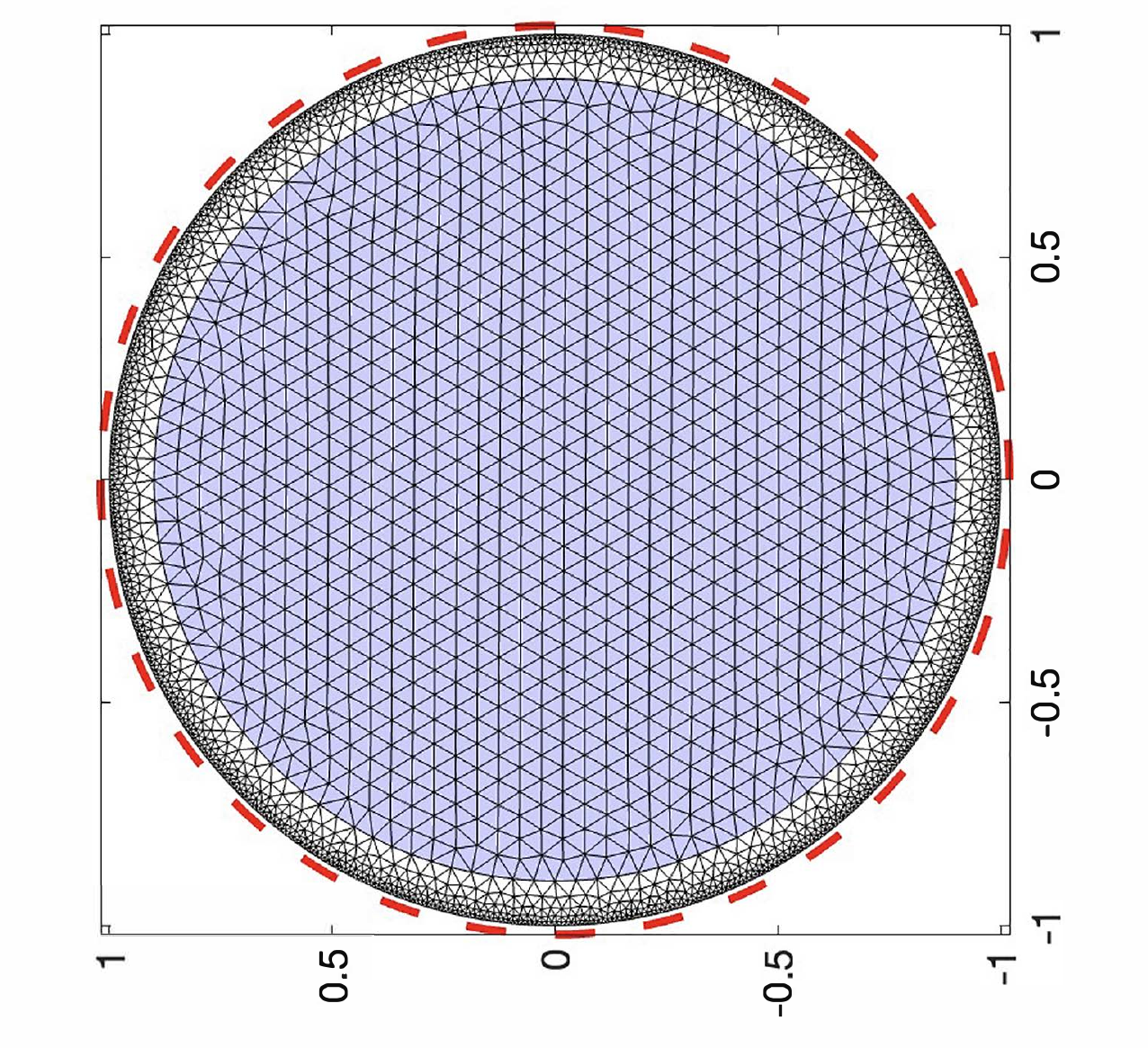}
\includegraphics[width=6.2cm,angle=-90]{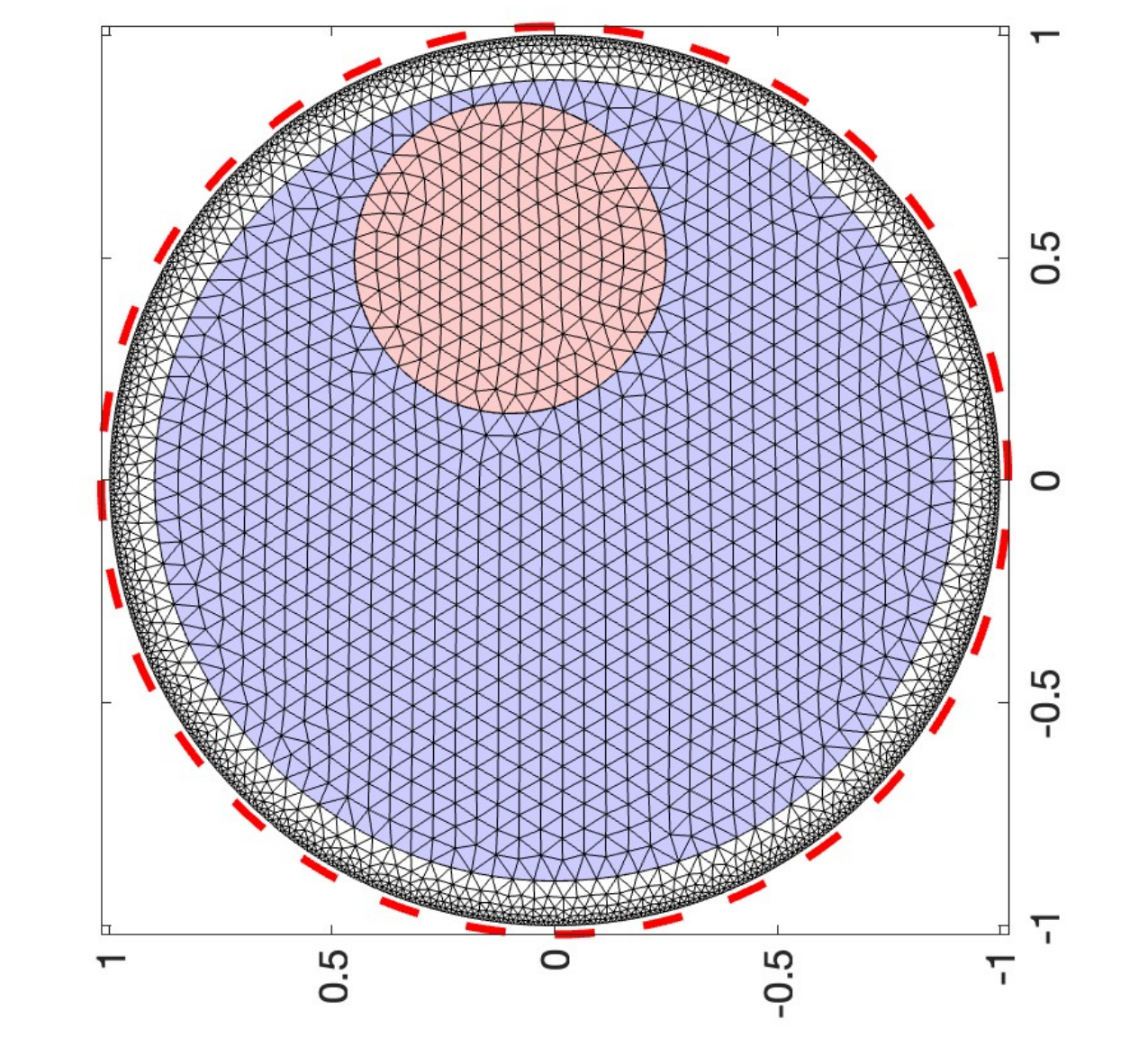}
\includegraphics[width=6.2cm,angle=-90]{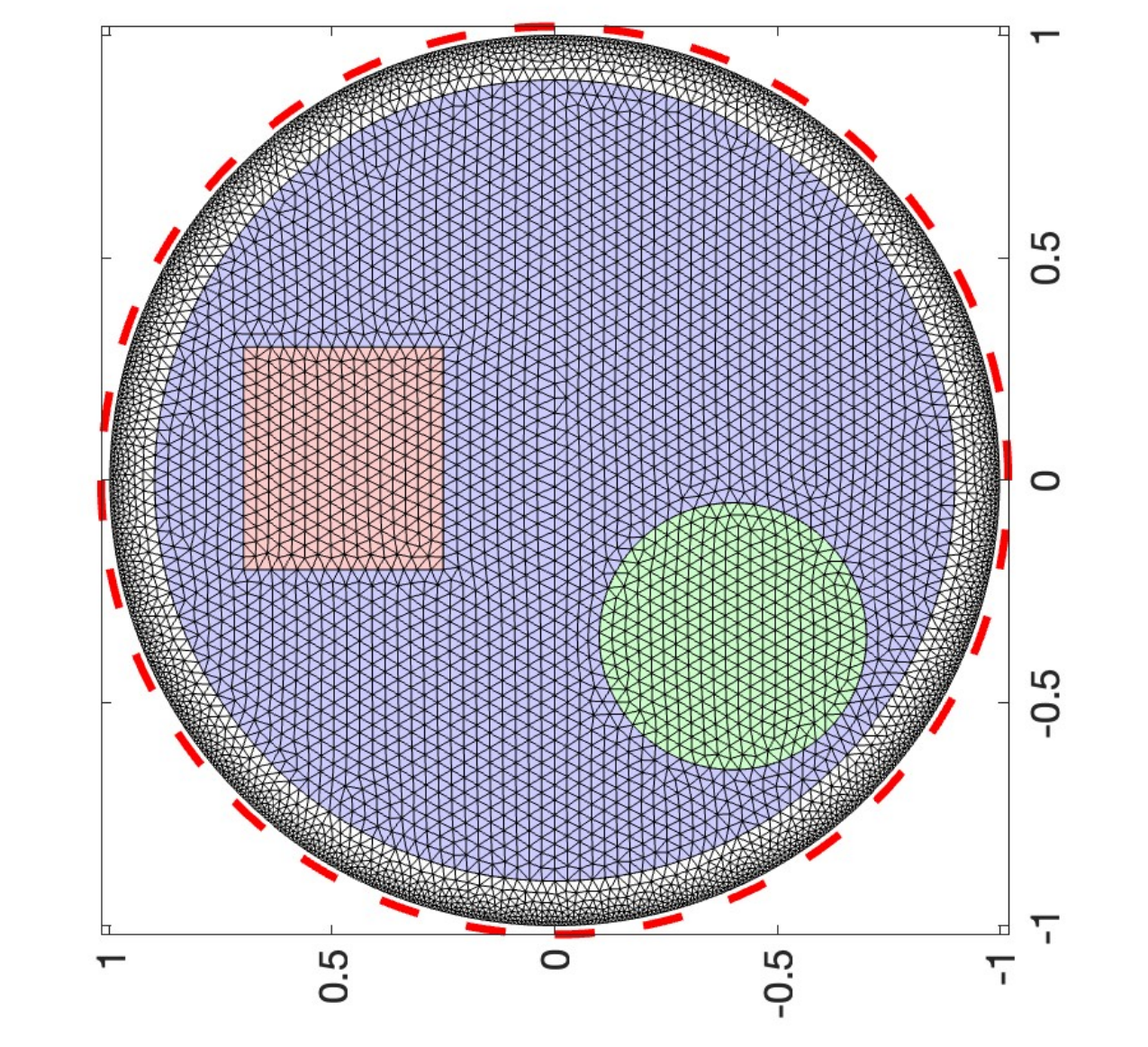}
}
\caption{ On the left, the tessellation used for solving the inverse problem is shown. The disc with blue background corresponds to the compact domain $D$ containing the support of $\delta\sigma$, while in the collar domain with no shading, the conductivity is assumed to be equal to a known constant $\sigma_0 = 1$. To capture the singularities at the edges of the electrodes indicated by red, a significant refining of the mesh near the boundary is necessary: The total number of elements in the mesh is 5\,816, while the number of elements in $D$, determining the number of degrees of freedom in the problem, is 1\,940. in the middle and on the right, the meshes used to generate the data. In the middle, the background value is equal to $\sigma_0 = 1$, while in the pink inclusion, the conductivity is set to $\sigma = 4.2$. On the right, the background conductivity is $\sigma_0=1$, and in the rectangular inclusion, $\sigma = 3.5$, while in the circular inclusion, the value is $\sigma = 0.4$. In the mesh on the right, the meshing near the electrodes is finer than in the mesh on the left to ensure that the reconstruction results are not mesh-dependent.}
\label{fig:meshes}
\end{figure}

The contact impedances of the electrodes are assumed to be all equal, $z_\ell = z_0$ and the value in the first numerical test is assumed to be known. In our simulations, we use the value $z_0 = 1\times 10^{-6}$.

In this simulation, the finite element mesh to be used for the reconstruction, shown in Figure~\ref{fig:meshes} (left), consists of 5\,816 elements, of which 1\,940 are inside the domain $D$, defining the dimension of the variable $\xi$. The number of edges defining the dimensionality of the unknown vector of increments is $N = 2\,970$.  In the simulations, the basis $\{\mathcal{E}_\ell\}_{\ell=1}^{L-1}$ of $\R^L_0$ is the trigonometric basis,
\[
 \big({\mathcal{E}}_\ell\big)_j =c_\ell \cos\frac{2\pi}{L}\ell j, \quad 1\leq \ell\leq \frac L2,\quad  \big({\mathcal{E}}_\ell\big)_j =c_\ell \sin\frac{2\pi}{L}\big(\ell-\frac L 2\big) j, \quad \frac L2+1\leq \ell<L,
\] 
where the constants $c_\ell$ are normalizing constants that make the basis orthonormal. The data correspond to a full frame, hence $b\in\R^m$ with $m=L(L-1) = 992$, roughly one third of the number $N$ of degrees of freedom. We generate two sets of simulated data, using the conductivities and FEM discretizations shown in Figure~\ref{fig:meshes} (center and right).

In the first numerical test, we use the simulated data generated by the model with one circular inclusion (center of Figure~\ref{fig:meshes}), corrupting the computed noiseless data by additive Gaussian scaled white noise with standard deviation $\omega$ equal to $0.1\%$ of the maximum of the noiseless voltages at the electrodes.
We run the IAS algorithm using the gamma hyperprior ($r=1$). The hyperparameter $\eta$ is set at
$\eta = 1\times 10^{-5}$ to guarantee strong sparsity promotion. To define the vector $\vartheta$, we compute the sensitivities at $\sigma = \sigma_0$, or $\xi = 0$, and define
$\vartheta$ by the formula (\ref{vartheta}), setting the constant $C$ so that $\max(\vartheta_j) = 4$, i.e., the expectations of the variancee $\theta_j$ are bounded by $\beta\max(\vartheta_j) \approx 6$.
We stop the IAS iterations when the condition (\ref{stopping}) with $\Delta_{\rm rel}\theta^\ell < \delta = 2\times 10^{-2}$ is satisfied. The condition is reached in sixteen iterations. Inside the non-linear least squares problem of updating $\zeta\in\R^N$, we perform only two linearization steps.  Numerical tests indicate that further linearizations do not improve significantly the results, not justifying the increase in the computing time for recalculating the Jacobian. With the current mesh size, the re-evaluation of the Jacobian using the adjoint formula requires approximately one second on an Apple M2 Ultra processor. 

\begin{figure}[ht!]
\centerline{
\includegraphics[width=17cm]{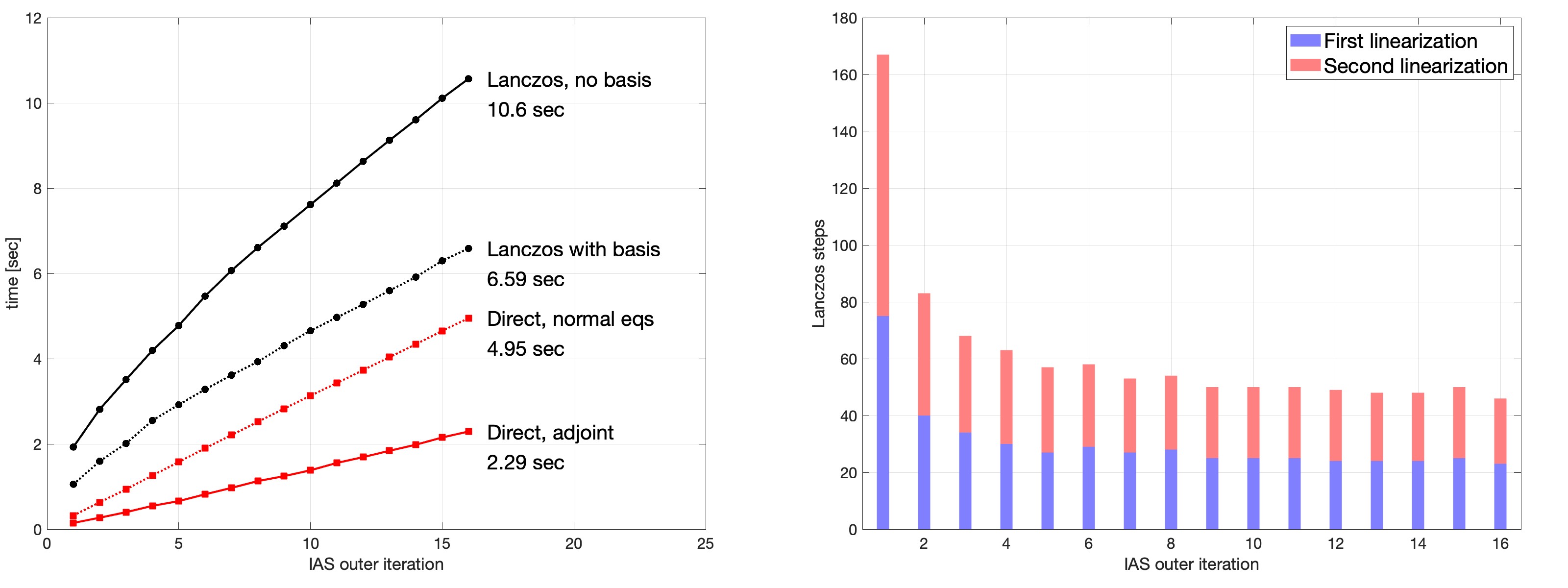}
}
\caption{ Left panel shows the cumulative times needed to solve the linearized least squares problems of the IAS algorithm. The times do not include the computation of the Jacobian, which takes about one second per re-evaluation. Here, ``Lanczos, no basis" refers to the Lanczos bidiagonalization without keeping the basis of the Krylov subspace in memory, while ``Lanczos with basis" saves the basis. ``Direct, normal eqs" refers to a direct solver of the system (\ref{normal}), and ``Direct, adjoint" to the direct solver of the adjoint system (\ref{adjoint}). 
The plot shows that for midsize problems in which the matrices are available and can be kept in the cache, the direct solvers are competitive, the adjoint formulation being the fastest. The calculations were performed using Apple M2 Ultra processor in a Mac Studio computer. The right panel shows the number of Lanczos steps needed to solve the linear system (\ref{adjoint}) at required accuracy of having the norm of the residual below the threshold of $10^{-8}$. }
\label{fig:timings}
\end{figure}

We test the computing time of different solvers of the linearized problem, comparing the performance of the direct solver of the normal equations (\ref{normal}), direct solver of the adjoint problem (\ref{adjoint}) in the data space, and Lanczos bidiagonalization process applied to the adjoint problem. Moreover, we compare two versions of the Lanczos process, one which keeps the orthonormal basis of the Krylov subspace in memory, and one which computes the solution in flight without keeping the vectors in memory. The latter one is bound to be slower as it has to run the process twice, however, being more parsimonious in memory requirement, it can be  used in computing environments where limited memory is an issue. Figure~\ref{fig:timings} shows the cumulative computing times needed to perform the linear solutions. As is evident from the results, the direct solvers are performing better here as the dimensions of the matrices are small enough to be kept in the cache, however, there is a significant speedup in using the adjoint problem formulation.  The right panel of Figure~\ref{fig:timings} shows the number of Lanczos iterations needed to solve the linear system with prescribed accuracy, set here at norm of the residual error smaller than $10^{-8}$. We observe that as the solution becomes sparser, the effective dimension of the Krylov subspace drops   to 22 or 23. 

\begin{figure}[ht!]
\centerline{
\includegraphics[width=\textwidth]{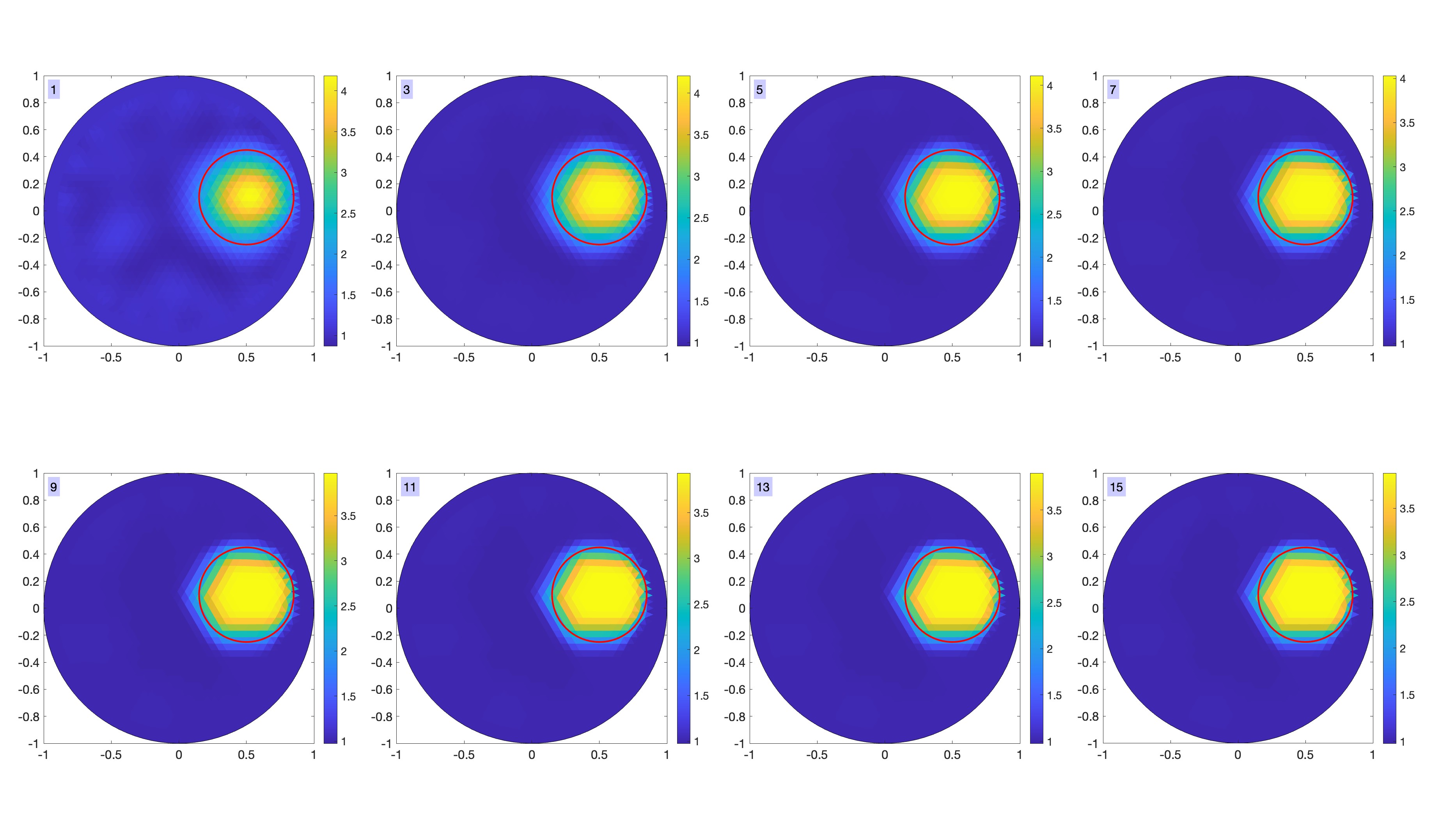}
}
\caption{ The iterative solutions estimating the conductivity $\sigma$ after an odd number of IAS iterations with the hyperparameter model with $r=1$ corresponding to the gamma distribution prior model. Each non-linear least squares solution involves two linearization steps. Observe that after the eighth iteration (lower row), the estimated conductivity remains practically unaltered. 
}
\label{fig:conductivities}
\end{figure}

The reconstructions shown in Figure~\ref{fig:conductivities} indicate that while the variance parameter $\theta$ keeps slightly changing, the effect on the reconstructions of the conductivity after eight iterations is minimal. In fact, the sharpness of the image does not improve significantly, but the dynamical range seems even to deteriorate slightly. A significant feature of the reconstructions is that the reconstructed inclusion is conformal with the underlying mesh. This is not surprising, as the ill-posedness of the EIT problem makes the likelihood insensitive to small shape variations, in particular further away from the electrodes, and therefore the grid-bound prior effectively determines the shape of the inclusion.

\begin{figure}[ht!]
\centerline{
\includegraphics[width=17cm]{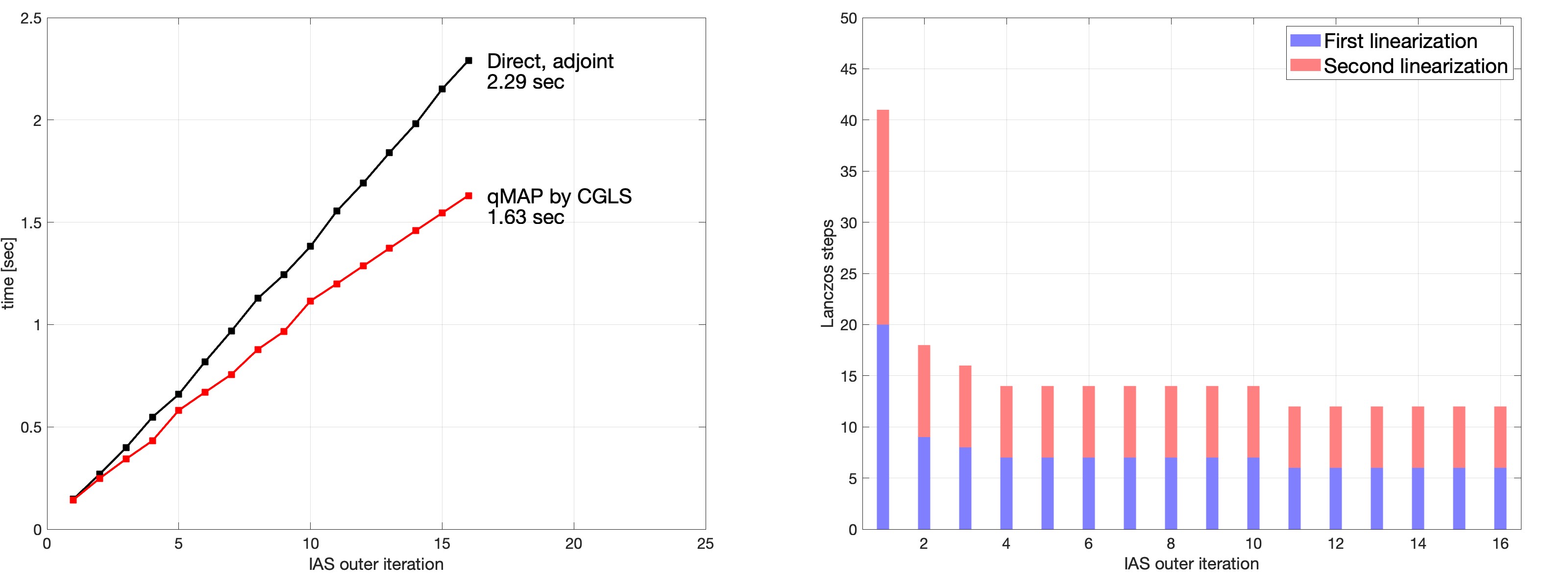}
}
\caption{The left panel shows the cumulative times needed to solve the linearized least squares problems of the IAS algorithm using the direct method of the adjoint map, compared to the approximate CGLS solution with early stopping at discrepancy.  The right panel shows the number of CGLS steps (equivalent to the number of Lanczos steps) needed before the stopping condition is met.}
\label{fig:timings MAP vs qMAP}  
\end{figure}

\begin{figure}[ht!]
\centerline{
\includegraphics[height=5cm]{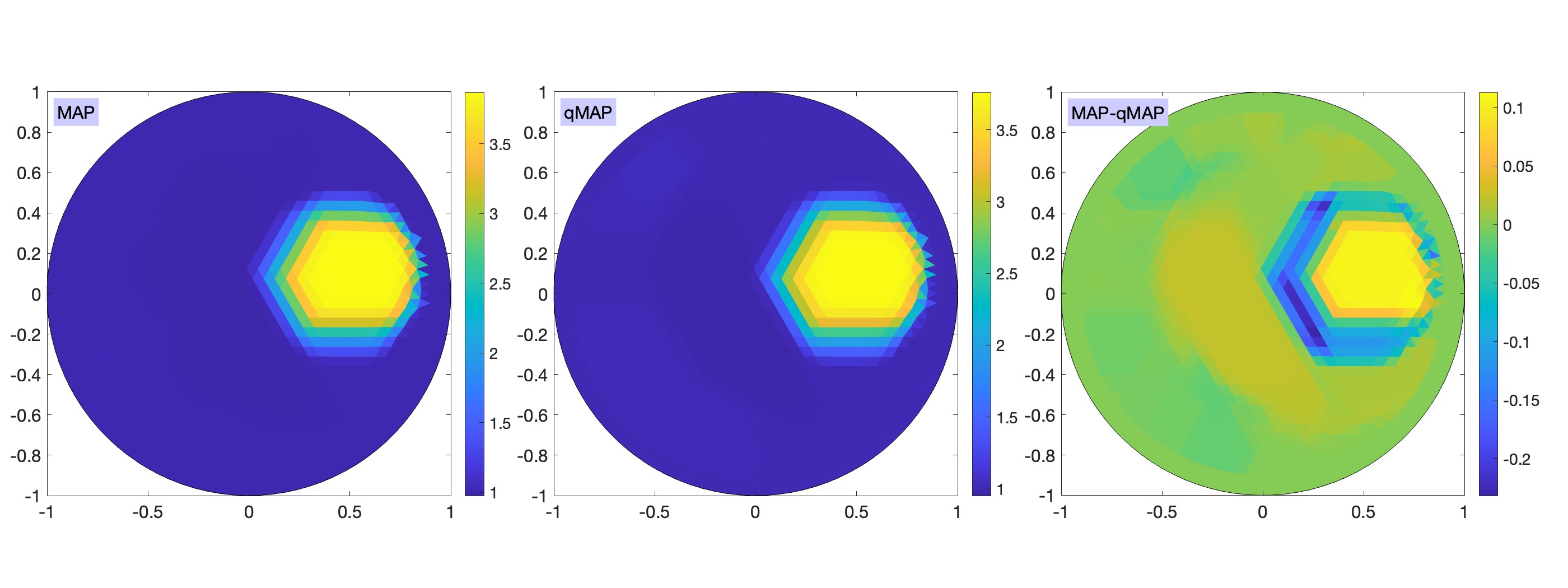}
}
\caption{ Comparison of the MAP estimate (left) and the qMAP estimate (center) obtained by replacing the exact solution of the approximation replacing the Tikhonov-regularized least squares problem by the least squares problem solved by using CGLS with an early stopping at the discrepancy. On the right, the difference of the two estimates is shown.}
\label{fig:sigma MAP vs qMAP}
\end{figure}

While solving the MAP estimate in the current simulation is numerically viable, it is interest to consider the performance of the qMAP approximation in which the linear leat squares problem with Tikhonov-type regularization is replaced by the approximation by CGLS with an early stopping rule. Figure~\ref{fig:timings MAP vs qMAP} shows the comparison between the computing times using the fastest solver for the exact problem, the direct solver of the adjoint problem, versus those needed for the CGLS approximation. The comparison shows that the qMAP is slightly faster to solve, and the number of iterations given in the same figure needed for reaching the convergence is significantly lower than the number of Lanczos steps in the computation of the MAP estimate. Moreover, we compare the accuracy of the qMAP approximation: Figure~\ref{fig:sigma MAP vs qMAP} shows the final estimates after 16 outer IAS iterations, along with the difference of the estimates. The results show that visually, the reconstruction images are essentially identical, however, the difference image shows that the qMAP estimate has a more limited dynamical range, roughly 3\%, with a faint ringing pattern around the inclusion.

\begin{figure}[ht!]
\centerline{
\includegraphics[height=5cm]{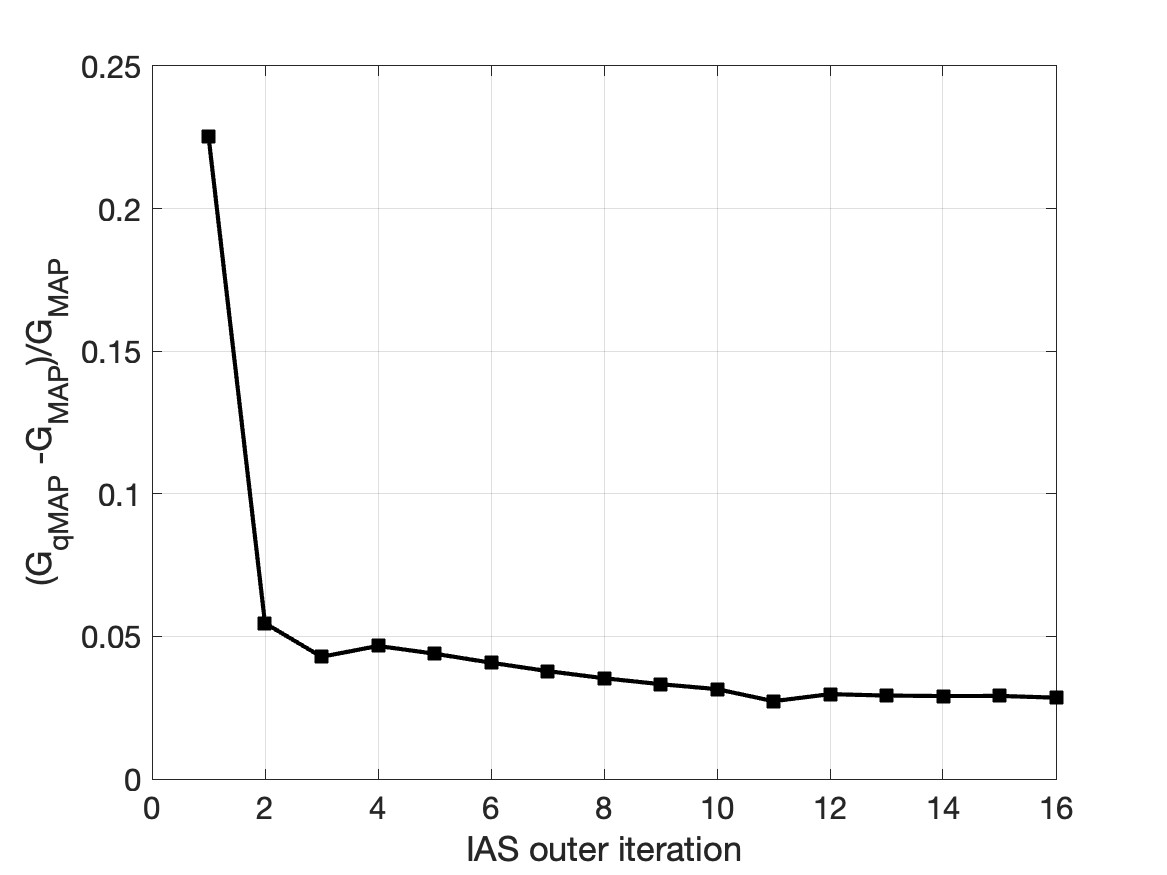}
}
\caption{ The relative difference of the Gibbs energy to be minimized when the solution of the Tikhonov regularized linear problem is replaced by the approximate solution using CGLS with an early stopping.}
\label{fig:Gibbs}
\end{figure}

Finally, Figure~\ref{fig:Gibbs} compares the values of the Gibbs energy objective function (\ref{Gibbs}) evaluated during the IAS iterations for the MAP and qMAP estimates. 
Denoting by $G_{\rm MAP}^\ell$ and $G_{\rm qMAP}^\ell$ the values of the Gibbs energy at the $\ell$th iteration while computing the approximations of the MAP and qMAP estimate, respectively, we expect that $G_{\rm MAP}^\ell \leq G_{\rm qMAP}^\ell$.  In the figure, we plot the values of the relative differences,
\[
 \Delta G^\ell = \frac{G_{\rm qMAP}^\ell - G_{\rm MAP}^\ell}{G_{\rm MAP}^\ell}, \quad ,\ell = 1,2,\ldots.
\] 
After a few iterations, the relative difference drops under 5\%, reaching asymptotically a level of approximately 3\%, in accordance with the expectation.

\begin{figure}[ht!]
\centerline{
\includegraphics[width=\textwidth]{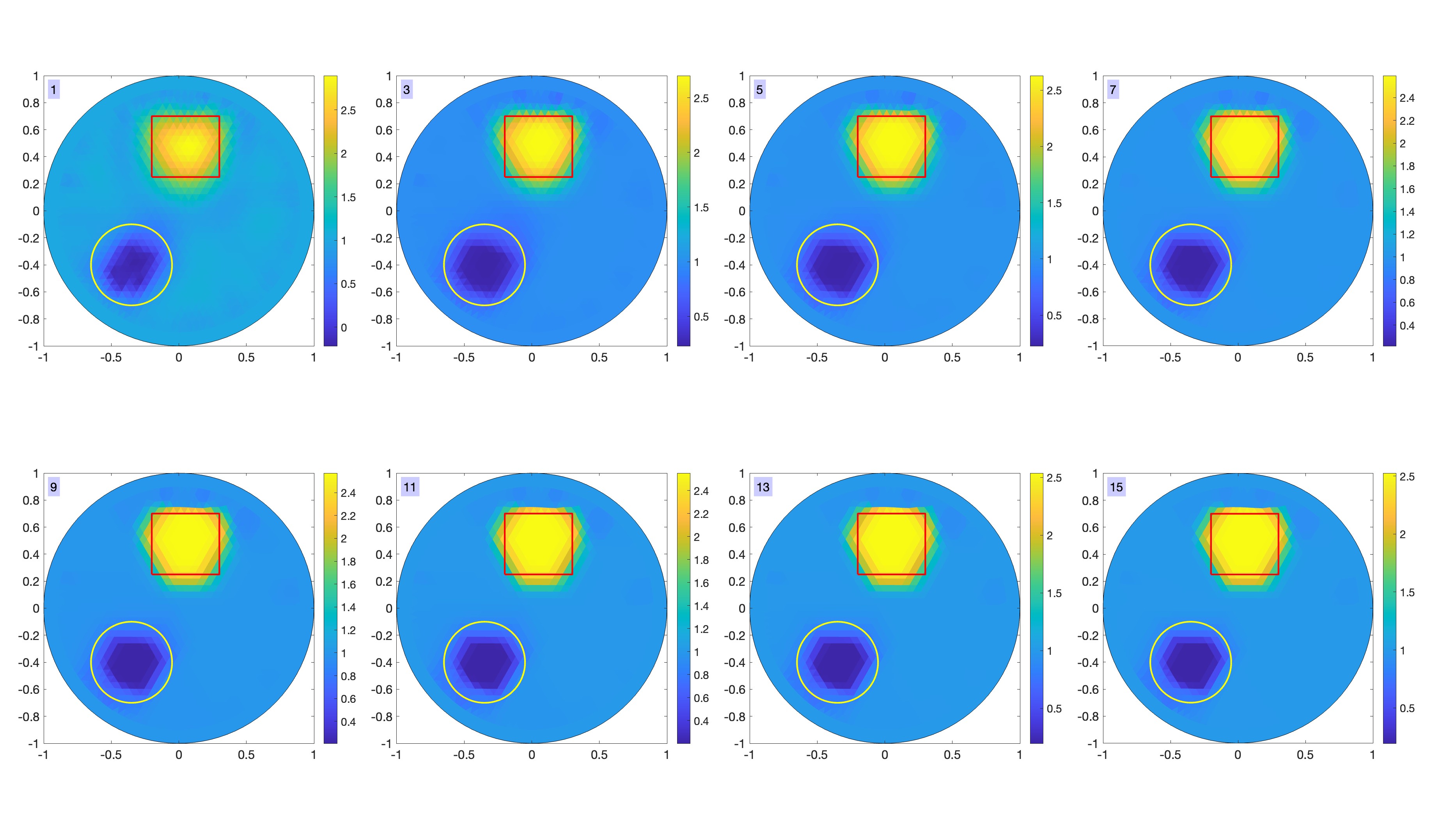}
}
\caption{The iterative solutions of the conductivities after an odd number of IAS iterations. Each non-linear least squares solution involves two linearization steps. Observe that after the eighth iteration (lower row) the estimated conductivity remains practically unaltered. }
\label{fig:conductivities2}
\end{figure}

In the second experiment, we run the IAS algorithm using the data generated by the generative conductivity of two inclusions shown on the right of Figure~\ref{fig:meshes}.
The noise level and model parameters are kept equal to the ones in the first experiment. Using the gamma hyperprior model, the algorithm converges again in 16 iterations: A few selected approximations are shown in Figure~\ref{fig:conductivities2}. The results are qualitatively similar to the ones obtained in the first experiment: The shape of the inclusions, in particular along the edges far away from the electrodes, is primarily determined by the underlying mesh. Numerical experiments using a finer reconstruction grid (not shown here) indicate that the phenomenon is not due to use of a too coarse mesh.

\begin{figure}[ht!]
\centerline{
\includegraphics[width=\textwidth]{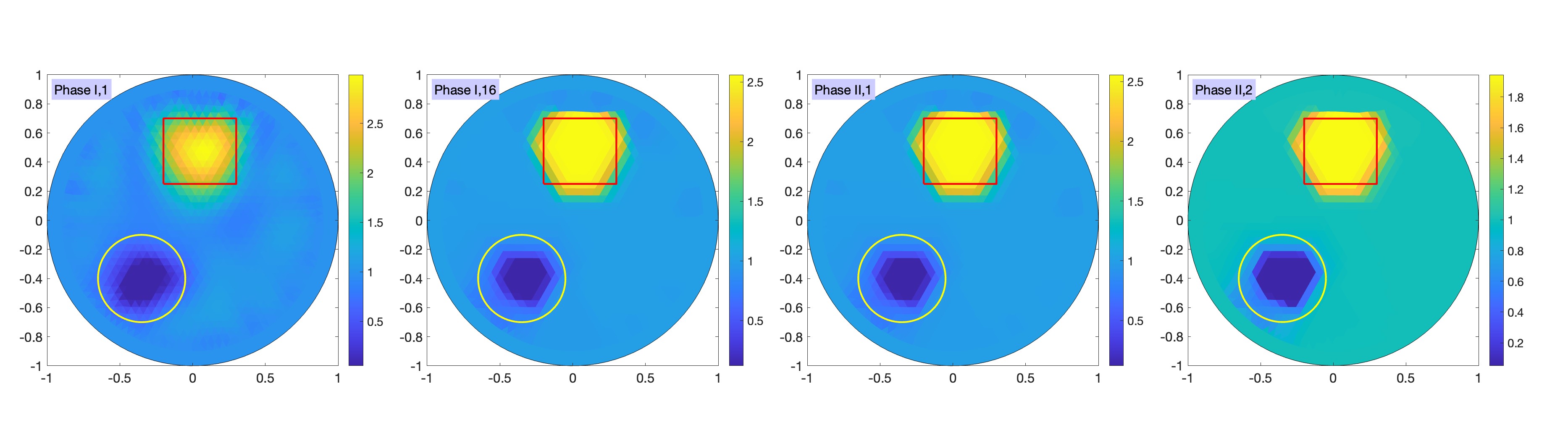}
}
\caption{ Estimates of the conductivities using the hybrid IAS scheme. From left to right: First iterate using the gamma hypermodel ($r=1$); iteration number 16 after the algorithm has reached convergence with the gamma hypermodel; First update after the switch to the generalized gamma hypermodel ($r=1/2$, referred her to as Phase II); and the second update with the generalized model.}
\label{fig:conductivities3}
\end{figure}

Finally, we test the effect of the hybrid scheme by switching the prior model from the case $r=1$ gamma distribution to a greedier generalized gamma model, setting $r = 1/2$.  Figure~\ref{fig:conductivities3} shows the two first iterations using the generalized gamma model. The greedier model clearly flattens the estimate in the background and within the inclusions, however the greedier model emphasizes even more the underlying grid structure and deteriorates the dynamic range. Running the iterations to convergence (not shown) deteriorates the dynamic range, possibly indicating that the algorithm converges to a local minimum.

\section{Conclusions}

Sparsity promoting hierarchical priors for identifying sparse blocky targets in inverse problems have previously shown to lead to numerically efficient computational algorithms for linear inverse problems. Their potential for solving non-linear problems was demonstrated in the article \cite{manninen2022sparsity} by an application to diffuse optical tomography. In this article, the methodology is developed further in the context of electrical impedance tomography, with a particular emphasis on the performance of different alternatives for solving the embedded linear least squares problem of Tikhonov type that is a subtask in the iterative IAS process. Moreover, the performance of the quasi-MAP estimate previously advocated in  articles of the authors was assessed in the context of the EIT. The numerical tests indicate that compared to the MAP estimate, the quasi-MAP estimate is computationally more lightweight, and while qualitatively very similar to the MAP estimate, the performance in terms of the dynamical range is slightly inferior.  The present work also addresses the question of convergence of the IAS algorithm in the current non-linear setting, and while a definitive uniqueness result could not be established, the analysis gives a partial convexity estimate of the Gibbs energy. Further studies to establishing whether it is possible to find conditions on model that  guarantee global uniqueness will be the topic of future studies. 

\section*{Acknowledgements}

The authors acknowledge the partial support by the NSF, grants DMS 1951446 for Daniela Calvetti and DMS 2204618 for Erkki Somersalo. Furthemore, Erkki Somersalo is grateful for the support of John Simon Guggenheim Foundation. 
Monica Pragliola acknowledges the National Group for Scientific Computation (INdAM-GNCS), Research Projects 2023, and the FRA (Fondi Ricerca Ateneo, University of Naples Federico II) project `HyRED'.

\bibliographystyle{siam} 
\bibliography{biblio_EIT}
\end{document}